\documentclass{imsart}
\RequirePackage{natbib}
\usepackage[english]{babel}
\usepackage[latin1]{inputenc}
\usepackage[T1]{fontenc}

\usepackage{graphicx}
\usepackage{pgf,pgfnodes,pgfarrows,pgfshade}
\usepackage{graphics}
\usepackage{caption}
\usepackage{amsmath}
\usepackage{amssymb}
\usepackage{amsthm}
\usepackage{hyperref}
\usepackage{array}
\usepackage{color}

\startlocaldefs

\theoremstyle{plain}
\newtheorem{theo}{Theorem}[section]
\newtheorem{lem}[theo]{Lemma}
\newtheorem{prop}[theo]{Proposition}
\newtheorem{corollary}[theo]{Corollary}
\newtheorem{remark}{Remark}

\theoremstyle{definition}

\theoremstyle{remark}

\newcommand{\E}{\mathbb{E}}
\newcommand{\p}{\mathbb{P}}
\newcommand{\R}{\mathbb{R}}

\newcommand{\N}{\mathbb{N}}

\endlocaldefs

\begin{document}

\title{A statistical analysis of a deformation model with Wasserstein barycenters : estimation procedure and goodness of fit test.}
\runtitle{Deformation model for distributions with Wasserstein barycenters}

\begin{frontmatter}
\begin{aug}
\author{\fnms{Del Barrio} \snm{Eustasio}\ead[label=e1]{tasio@eio.uva.es}}, \and
\author{\fnms{Lescornel} \snm{H\'el\`ene}\ead[label=e2]{helene.lescornel@math.univ-toulouse.fr}}, \and
\author{\fnms{Loubes} \snm{Jean-Michel}\ead[label=e3]{loubes@math.univ-toulouse.fr}}
%

\affiliation{Institut de Math\'ematiques de Toulouse}

\address{Address of the First and Second authors\\
IMT \\
118 route de Narbonne \\
31000 Toulouse \\
FRANCE\\
\printead{e1}\\
\phantom{E-mail:\ }\printead*{e3}}

\runauthor{E. del Barrio, H. Lescornel and J-M. Loubes}

\end{aug}

\begin{abstract}{ }
 We propose a study of a distribution registration model for general deformation functions. 
In this framework, we provide estimators of the deformations as well as a goodness of fit 
test of the model. For this, we consider a criterion which studies the Fr\'echet mean 
(or barycenter) of the warped distributions whose study enables to make inference on the model. 
In particular we obtain the asymptotic distribution and a bootstrap procedure for the  Wasserstein variation.
\end{abstract}


%

\end{frontmatter}


\section{Introduction}\label{s:intro}
Giving a sense to the  notion of {\it mean behaviour} may be counted among the very early activities of statisticians. 
When confronted to large data sample,  the usual notion of Euclidean mean is too rough since the information conveyed 
by the data possesses an inner geometry far from the Euclidean one. Indeed, deformations  on the data such as 
translations, scale location models for instance or more general warping procedures prevent the use of the 
usual methods in data analysis. This problem arises naturally for a wide range of statistical research fields 
such as functional data analysis for instance   \cite{Ramsay-Silverman-05}, \cite{bercufraysseaos}  and 
references therein, image analysis in \cite{TrouveY05}  or \cite{JASA}, shape analysis in \cite{kendall}, \cite{bb206} 
or \cite{MR2640651}, with many applications ranging from biology in \cite{Bolstad-03} to pattern recognition \cite{Sakoe-Chiba-78} just to name a few.

\indent The same kind of issues arises when considering the estimation of distribution functions observed with deformations.  
This situation occurs often in biology,  for example when considering gene expression. However, 
when dealing with the registration of warped distributions, the literature is scarce. We mention 
here the method provided for biological computational issues known as quantile normalization in~\cite{Bolstad-03} and 
the related work \cite{GALLON-2011-593476}. Very recently using optimal transport methodologies,  comparisons of 
distributions have been studied using a notion of Fr\'echet mean for distributions, see for instance in 
\cite{agueh2010barycenters} or a notion of depth as in \cite{2014arXiv1412.8434C}. 

Actually the observations are said to come from a deformation model if they can be written 
as $$ X_{i,j}=  \left(\varphi^\star_j\right) ^{-1}  \left(\varepsilon_{i,j}\right),$$  for 
$j=1,\dots,J$ where $\left(\varepsilon_{i,j}\right)$ defined for all ${ 1 \leqslant i \leqslant n, 1 \leqslant j \leqslant J} $ 
are  i.i.d. random variables with unknown distribution $\mu$ and for deformation functions $\varphi_j^\star$. 
This model is the natural extension of the functional deformation models studied in the statistical literature 
for which estimation procedures are provided in \cite{Gamboa-Loubes-Maza-07} while testing issues are tackled 
in \cite{MR3323102}. Within this framework, statistical inference on deformation models for distributions have been studied first 
in~\cite{MR1704844}, ~\cite{MR1625620} and \cite{Freitag2005123}, where tests are provided in the case of parametric functions, while 
the estimation of the parameters is studied in \cite{agullo2015parametric}. 

In this work, after recalling the model we use in Section~\ref{s:model}, we tackle the problem of providing a goodness of fit test in a general non parametric deformation model. For this, we will use an alignment criterion with respect to the Wasserstein's barycenter of a deformation of the observed distributions. This requires an equivalent of a central limit theorem for the Wasserstein variation of a barycenter of measures in both the general case in Section~\ref{s:nonparam} and under the null assumption (observations are drawn from the deformation model)  in Section~\ref{s:param}. We obtain the asymptotic distribution of 
the matching criterion in both cases, with a different normalization under the null assumption (only for the parametric case). 
For this, we will need to build estimates of the deformation parameters with respect to this particular alignment criterion and study their 
behavior in Section~\ref{s:s:estim}. Finally testing procedures are given in Section~\ref{s:test}. They rely on 
the estimation of the quantiles of the empirical process of the Wasserstein's variation which is obtained using a bootstrap procedure proved 
in Section~\ref{s:boot}. Proofs are postponed to  Section~\ref{s:append}.

\section{A deformation model for distributions} \label{s:model}
Assume we observe $J$  samples of $n$ i.i.d random variables $X_{i,j}$ with 
distribution $\mu_j$,  associated to a distribution function $F_j: ( c_j, d_j)  \mapsto (0,1)$ with density with respect to the Lebesgue measure $f_j$. 
Let $\mu_{n,j}$ and $F_{n,j}$ be the empirical measure and empirical distribution function associated to the sample 
$\left( X_{i,j} \right)_{ 1 \leqslant i \leqslant n}$.

Our aim is to test the existence of a distribution's deformation model, in the sense that  all  
the distributions $\mu_j$ would be warped from an unknown distribution template $\mu$  by a  deformation 
function $\varphi^\star_j$. More precisely, consider a family of  warping functions $\mathcal{G}=\mathcal{G}_1\times \dots \times\mathcal{G}_J $  such that
\begin{align}
 \label{hyp:fct} \tag{$\mathbf{\mathcal{A}1}$} \text{ For all } h \in \mathcal{G}_j, 
h : \begin{array}{c}  (c_j , d_j )\rightarrow (a , b )\\ x \mapsto h \left( x\right)
\end{array} \text{ is invertible, increasing,} \\ \notag \text{ and s.t.}
-\infty \leqslant a  < b \leqslant +\infty, \quad -\infty \leqslant c \leqslant c_j < d_j\leqslant d \leqslant +\infty.
\end{align}
We would like to build a goodness-of-fit testing procedure for the following model
\begin{align}
\label{eq:H0}
\text{ There exist }& \left(\varphi^ \star_1, \dots \varphi^\star_J \right)\in \mathcal{G} \text{ and } \left(\varepsilon_{i,j}\right)_{\substack{1 \leqslant i \leqslant n \\ 1 \leqslant j\leqslant J}} \text{ i.i.d. such that } \notag \\
& X_{i,j}=\left(\varphi^\star_j \right)^{-1} \left(\varepsilon_{i,j}\right)\quad \forall 1 \leqslant j \leqslant J\tag{$\mathcal{H}$}
\end{align} 
Denote by  $G$ the distribution function of $\varepsilon$ with law  $\mu$ with support $(a,b)$, while  $G_{n,j}$  is the corresponding  empirical version.
 
Our criterion will be based on the Wasserstein distance $W^2_2$ 
since this distance is well suited to compared deformations between distributions. For $d\geqslant 1$,  consider  the following set 
$$\mathcal{W}_2\left( \R^d\right) =\left\{ P \text{ probability on }\R^d \text{  with finite second moment}\right\}.$$

For two probabilities $\mu$ and $\nu$ in $\mathcal{W}_2\left( \R^d \right)$ , we denote by $\Pi(\mu, \nu)$ the set of all
probability measures $\pi$ over the product set $\R^d \times\R^d $
with first (resp. second) marginal $\mu$ (resp. $\nu$).

The transportation cost with quadratic cost function, or quadratic transportation cost,
between these two measures
$\mu$ and $\nu$   is defined as

\begin{equation*}
\label{eq:infwasser}
 \mathcal{T}_2(\mu, \nu) = \inf_{ \pi \in \Pi(\mu, \nu)} \int \left\Vert x - y\right\Vert ^2 d \pi(x,y).
\end{equation*}
The quadratic transportation cost allows to endow the set $\mathcal{W}_2\left(\R^d\right)$
 with a metric by defining  the $2$-Wasserstein distance between $\mu$ and $\nu$ as
$ W_2(\mu, \nu)  = \mathcal{T}_2(\mu, \nu) ^{1/2}.$
More details on Wasserstein distances and their links with optimal transport problems can be found in  
\cite{rachev} or \cite{villani2009optimal} for instance.

%
%
%
%

Here we will consider probabilities in $\mathcal{W}_2\left( \R\right)$. In 
this case,  the Wasserstein distance can be written as
\begin{equation}
\label{distW}
W_2^2\left( \mu, \nu  \right) = \int_0^1 \left( F^{-1} \left( t \right) -G^{-1} \left( t \right) \right)^2 dt,
\end{equation}
where $F$ (resp. $G$) is the distribution function associated with $\mu$ (resp. $\nu $).  

Moreover, we are  dealing with more than two probabilities and so we are interested in a global measure
of separation. So consider   the  \textit{Wasserstein $2$-variation} of $\nu_1,\ldots,\nu_J$,  
defined as follows. Given probabilities $\nu_1,\ldots,\nu_J$ on $\mathbb{R}^d$ with finite $2$-th moment let
$$V\left( \nu_1, \dots, \nu_J \right) = \inf_{\eta\in \mathcal{W}_2\left(\mathbb{R}^d\right)} \left(\frac 1 J \sum_{j=1}^J {W}_2^2 (\nu_j,\eta) \right)^{1/2}$$
 be the Wasserstein $2$-variation
of $\nu_1,\ldots,\nu_J$. In~ \cite{agueh2010barycenters}, the minimizer  of $\eta \mapsto \frac 1 J \sum_{j=1}^J {W}_2^2 (\nu_j,\eta) $ is proved to exist. This measure $\nu_B $ is called the barycenter or Fr\'echet mean of $\nu_1,\ldots,\nu_J$. The authors prove properties of existence and uniqueness for  barycentres of measures in $\mathcal{W}_2\left(\mathbb{R}^d\right)$, while the properties of the empirical version are provided in~\cite{boissard2014}.
 
We propose to use this Wasserstein $2$-variation as a goodness of fit criterion for model~\eqref{eq:H0}. 
Since the true distribution $\mu$ is unknown, we first try to invert the warping operator and thus compute 
for each observation its image through a candidate  deformation $\varphi_j$, $$ Z_{i,j}\left(\varphi_j\right)= 
{\varphi_j} \left( X_{i,j} \right) \quad 1 \leqslant i \leqslant n,\quad 1 \leqslant j \leqslant J.$$ 
Note that $  Z_{i,j}\left(\varphi_j\right) \sim \mu_j(\varphi_j)$ with distribution function 
(under Assumption \ref{hyp:fct}) $F_j \circ \varphi^{-1}_j:=F_{\varphi_j}$. { Now, if we set $\varphi= 
\left( \varphi_1, \dots, \varphi_J \right) \in \mathcal{G}$, then the 
Fr\'echet mean of $\left(\mu_j(\varphi_j) \right)_{1 \leqslant j \leqslant  J}$ is the probability $\mu_B(\varphi)$
with quantile function}
$$F^{-1}_B(\varphi)\left( t \right)=\frac{1}{J}\sum_{k=1}^J \varphi_k\circ F_k^{-1}\left( t \right),$$  
{(see \cite{agueh2010barycenters}). 
We will write $\mu_{n,j}(\varphi_j)$ for the empirical measure on $Z_{i,j}(\varphi_j), 1\leq i\leq n$ and $\mu_{n,B}(\varphi)$ for the
corresponding Fr\'echet mean.}
It is important to remark that $$ {\rm under} \: \eqref{eq:H0} \quad \mu_B(\varphi^\star)=\mu= \mu_j(\varphi_j^\star), \: \forall 1 \leqslant j \leqslant J.$$
{Hence, } a natural idea to test whether \ref{eq:H0} holds, is to consider the Wasserstein 
$2$-variation of the $\left(\mu_j(\varphi_j)\right),\: {1 \leqslant j \leqslant J }$, that is to say the 
minimum alignment of the candidate warped distributions $\left(\mu_j(\varphi_j)\right)_{1 \leqslant j \leqslant J }$  
with respect to their { barycenter}, namely $\mu_B(\varphi)$. This optimization  program corresponds  
to the minimization { in $\varphi \in \mathcal{G}$} of the following theoretical criterion 
\begin{equation*}
U \left( \varphi \right):=V^2\left( \mu_1(\varphi_1), \dots, \mu_J(\varphi_J) \right)  = \frac{1}{J} \sum_{j=1}^J W_2^2(\mu_j(\varphi_j),\mu_B(\varphi)).
\end{equation*}
Its empirical version is given by {
$U_n\left( \varphi \right) = \frac 1 J W_2^2(\mu_{n,j}(\varphi_j), \mu_{n,B}(\varphi))$.
Inference about model \ref{eq:H0} can be based on the statistic $\inf_{\varphi \in \mathcal{G}} U_n\left(\varphi\right)$. 
In the next sections we analyse the behavior of this statistics under different setups.}


\section{Non parametric model for deformations} \label{s:nonparam}
{ We provide in the section a CLT for $\inf_{\varphi \in \mathcal{G}} U_n\left(\varphi\right)$ under
the following set of assumptions.}

\begin{equation}
\label{hyp:tassio13} \tag{$\mathbf{\mathcal{A}2}$} \mbox{For all $j$, $F_j$ is $C^2$ on  $( c_j; d_j)$, $f_j(x)>0$ if $x \in ( c_j; d_j)$ and}
\end{equation}
$$\sup_{c_j < x < d_j}\textstyle  \frac{F_j(x)\left( 1- F_j(x)\right) f'_j(x)}{f_j(x)^ 2} < \infty.$$
\begin{equation}
\label{hyp:Rajput}\tag{$\mathbf{\mathcal{A}3}$}
\text{ For some } q>1 \text{ and all } 1\leqslant j \leqslant J,  \quad
\textstyle\int_0 ^1 \frac{\left(t \left(1-t\right)\right)^\frac{q}{2}}{\left(f_j\left( F_j^{-1}\left( t \right) \right)\right) ^q}dt < + \infty
\end{equation}

For $q$ as in \ref{hyp:Rajput}, we set $p_0=\max\left(\frac{q}{q-1}, 2 \right) $ and define on 
$\mathcal{H}_j=C^1 (c_j;d_j) \cap  L^{p_0} \left( X_j \right) $ the norm
$\left\Vert h_j \right\Vert_{\mathcal{H}_j}  = \sup_{(c_j;d_j)} \vert h_j' \vert + 
\E \left[  \left\vert h_j\left(X_j \right) \right\vert ^{p_0} \right]^\frac{1}{p_0},$ and on 
the product space $\mathcal{H}_1 \times \dots \times\mathcal{H}_J$, $\left\Vert h \right\Vert_{\mathcal{H}}  
= \sum_{j=1}^J \left\Vert h_j \right\Vert_{\mathcal{H}_j}.$
The we make the following additional assumptions.

\begin{equation}
\label{hyp:regul}\tag{$\mathbf{\mathcal{A}4}$}
\mathcal{G}_j\subset \mathcal{H}_j \text{ is compact for } \Vert \cdot \Vert_{\mathcal{H}_j}
\text{ and } \sup_{h  \in \mathcal{G}_j}  \left\vert h^\prime (x^h _n)- h^\prime(x) 
\right\vert \underset{ \sup_{h  \in \mathcal{G}_j} \vert x_n^h- x  \vert \rightarrow 0}\rightarrow 0.
\end{equation}


\begin{equation}
\label{hyp:tassio12}\tag{$\mathbf{\mathcal{A}5}$}
 \text{for some }  r>4 \text{ and } 1 \leqslant j \leqslant J, \quad \E \left[ \left\vert X_j \right\vert^ r \right] < \infty 
\end{equation}
\begin{equation}
\label{hyp:tassio12theta}\tag{$\mathbf{\mathcal{A}6}$}
 \text{for some }  r>\mbox{max}(4,p_0) \text{ and }  1 \leqslant j \leqslant J, \quad \E \left[\sup_{h \in \mathcal{G}_j} \left\vert h \left( X_j\right)\right\vert ^ r \right]<\infty 
\end{equation}

Under \ref{hyp:fct} to \ref{hyp:tassio12theta}, we are able to provide the asymptotic distribution of $\inf_{\varphi \in \mathcal{G}}\sqrt{n} U_n(\varphi)$.
It is convenient at this point to give some explanation about the meaning of these assumptions. \ref{hyp:tassio13}
is a is a regularity condition on the distributions of the $X_j's$ (it holds, for instance, for Gaussian or Pareto distributions) required for strong
approximation of the quantile process, see \cite{csorgo83} for details. The
integrability condition \ref{hyp:Rajput}  is satisfied by the Gaussian distribution if $q<2$, see, e.g., \cite{rajput1972}.
\ref{hyp:regul} is related to the regularity of the deformation functions. Finally, \ref{hyp:tassio12}
and \ref{hyp:tassio12theta} are moment assumptions on the (possibly warped) observations.

\begin{theo}
\label{th:testgen}
Under  Assumptions \ref{hyp:fct} to \ref{hyp:tassio12theta}
\begin{eqnarray*}
\lefteqn{\sqrt{n} \Big( \inf_{\varphi \in \mathcal{G}}  U_n(\varphi) -\inf_{\varphi \in \mathcal{G}}U(\varphi) \Big)}\hspace*{1cm}\\
&\rightharpoonup &\inf_{\varphi  \in \Gamma}
\frac{2}{J} \sum_{j=1}^J \int_0^1  \varphi_j^\prime\circ F_j^{-1}\frac{B_j}{ {f_j\circ F_j^{-1}}} (\varphi_j \circ F_j^{-1} -F^{-1}_B(\varphi)),
\end{eqnarray*}
where $\Gamma=\{\varphi\in\mathcal{G}:\, U(\varphi)=\inf_{\phi \in \mathcal{G}}U(\phi)\}$ and 
$\left( B_j\right)_{1 \leqslant j \leqslant J}$ are independent Brownian bridges.
\end{theo}

A proof of Theorem \ref{th:testgen} is given in the Appendix below. We note that 
continuity of $U$ is follows easily from the choice of the norm on $\mathcal{G}$. 
Recall that $\mathcal{G }$ is compact and, consecuently, $\inf_{\varphi\in\mathcal{G}}U(\varphi)$
is attained. Hence, $\Gamma$ is a nonempty closed subset of $\mathcal{G}$ (in particular, it is also a compact set).
We note further that the random variables $\int_0^1  \varphi_j^\prime\circ F_j^{-1}\frac{B_j}{ {f_j\circ F_j^{-1}}} (\varphi_j \circ F_j^{-1} -F^{-1}_B(\varphi))$
are centered Gaussian, with covariance 
\begin{eqnarray*}
\lefteqn{\int_{[0,1]^2} (\min(s,t)- st)  {\textstyle\frac{ \varphi_{j}^{\prime}(F_j^{-1}(t)) }{ {f_j\left(F_j^{-1}(t) \right)}}} 
(\varphi_{j}( F_j^{-1}(t)) -F^{-1}_B(\varphi)(t))}\hspace*{2.5cm}\\
&\times &
{\textstyle\frac{ \varphi_{j}^{\prime}(F_j^{-1}(s)) }{ {f_j\left(F_j^{-1}(s) \right)}}} 
(\varphi_{j}( F_j^{-1}(s)) -F^{-1}_B(\varphi)(s))dsdt.
\end{eqnarray*}
In particular, if $U$ has a unique minimizer the limiting distribution in Theorem \ref{th:testgen} is normal. However,
our result works in more generality, even without uniqueness assumptions.

We remark also that although we have focused for simplicity on the case of samples of equal size,  
the case of different sample sizes, $n_j$, $ j =1,\dots,J$,
can also be handled with straightforward changes. If we assume
\begin{equation}\forall j :  n_j \rightarrow + \infty \text{ and } \frac{n_j}{n_1 + \dots + n_J} 
\rightarrow \left( \gamma_j \right) ^2>0,\label{eq:condi_ndif}\end{equation} 
then the result can be restated as 
\begin{align*}
\sqrt{\frac{n_1 \dots n_J}{\left(n_1 + \dots + n_J\right)^{J-1}}}  
\Big( \inf_{\mathcal{G}}U_{n_1,\dots,n_J} -
\inf_{\mathcal{G}}U\Big) \rightharpoonup  \inf_{\Gamma} \frac{2}{J} \sum_{j=1}^J \widetilde{S}_j,
\end{align*} 
where $U_{n_1,\dots,n_J}$ denotes the empirical Wasserstein variation computed from the samples and
$\widetilde{S}_j(\varphi)= \big(\Pi_{p\neq j} \gamma_p\big)\int_0^1  \varphi_j^\prime
\circ F_j^{-1}\frac{B_j}{ {f_j\circ F_j^{-1}}} (\varphi_j \circ F_j^{-1} -F^{-1}_B(\varphi)) $.

As a final remark in this section we note that 
in the case where \ref{eq:H0} holds, we have $\varphi_j\circ F_j^{-1}=F_B^{-1}(\varphi)$ for each $\varphi\in\Gamma$.
Thus, the result of Theorem \ref{th:testgen} becomes  $$\inf_{\varphi \in \mathcal{G}} \sqrt{n} U_n(\varphi) \rightharpoonup 0.$$
Hence, in this case we have to refine our study to understand well the  behavior of $\inf_{\mathcal{G}}  U_n$ when $n$ 
tends to infinity. This is what we consider in the next section. In this case we restrict
ourselves to a to a semiparametric warping model
where $\mu$ is unknown but where the deformations are indexed by a parametric family.

\section{A parametric model for deformations}\label{s:param}
In many cases, deformation functions can be made more specific in the sense that they follow a known shape 
depending on parameters that are different  for each sample. So consider the parametric model $\theta^\star=(\theta_1^\star,\dots,\theta_J^\star)$ such that 
$ \varphi^\star_j=\varphi_{\theta^\star_j},\quad {\rm for} \: \: {\rm all} \: \: j=1,\dots,J.$ 
Each $\theta^\star_j$ represents the warping effect that undergoes the $j^{\rm th}$ sample, which must be removed 
to recover the unknown distribution by inverting the warping operator. So Assumption \ref{eq:H0} becomes 
$$ X_{i,j}= \varphi_ {\theta^\star_j}^{-1} \left(\varepsilon_{i,j}\right),\: {\rm for} \: {\rm all} \: 1 \leqslant i \leqslant n,  1\leqslant j \leqslant J.$$
Hence, from now on, we will consider the following family of deformations, indexed by a parameter 
$ \lambda \in \Lambda \subset \R^p$: 
$$ \begin{array}{c c c} \varphi :  \Lambda \times (c;d) & \rightarrow & (a,b) \\
\left(\lambda,x\right) & \mapsto & \varphi_{\lambda} \left( x\right) \end{array} 
$$
Thus, the functions $U$ and $U_n$ are now defined on $\Theta= \Lambda^J$, and the criterion of interest 
becomes $\inf_{\lambda \in \Theta} U(\lambda)$. 
We also use the simplified notation $\mu_j(\theta_j)$ instead of $\mu_j\left( \varphi_{\theta_j} \right)$,
$F_B \left( {\theta} \right)$ for $F_B \left( \varphi_{\theta_1}, \dots, \varphi_{\theta_J} \right)$ and 
similarly for the empirical versions. Throughout this section \textit{ we assume that model
\ref{eq:H0} holds}. This means, in particular, that the d.f.'s of the samples, $F_j$,
satisfy $F_j=G\circ \varphi_{\theta_j^*}$, with $G$ the d.f. of the $\varepsilon_{i,j}$'s.

For the analysis of this setup, we adapt Assumptions  
\ref{hyp:fct} to \ref{hyp:tassio12theta}, replacing them by the following versions.
\begin{align}
 \label{hyp:fct_th} \tag{{\bf A1}} \text{ For all } \lambda \in \Lambda, 
\varphi_\lambda : \begin{array}{c} (c;d ) \rightarrow (a;b)\\ x \mapsto \varphi_{\lambda} \left( x\right)
\end{array} \text{ is invertible, increasing,} \\ \notag \text{ and s.t.}
-\infty \leqslant a  < b \leqslant +\infty, \quad -\infty \leqslant c \leqslant c_j < d_j\leqslant d \leqslant +\infty.
\end{align}
We replace \ref{hyp:tassio13} by: 
$G$ is $C^2$ with  $G'(x)=g(x)>0$  on $(a,b)$ and
\begin{equation}
\label{hyp:tassio13_eps} \tag{{\bf A2}}
\sup_{a< x <b} \frac{G(x)\left( 1- G(x)\right) g'(x)}{g(x)^ 2} < \infty
\end{equation}
Now, instead of \ref{hyp:Rajput} to \ref{hyp:tassio12} we assume
\begin{align}
\label{hyp:regulbis_th}\tag{{\bf A3}}
\varphi  &\text{  is continuous w.r.t. } x  \text{ and } \lambda  \\ \forall \lambda \in \Lambda\text{, }\varphi_\lambda   & \text{ is } C^ 1 \text{ with respect to } x\text{, } \Lambda  \text{ is compact}\notag
\end{align}
\begin{align}
\label{hyp:regul_th}
&d\varphi \text{ is bounded on } \Lambda \times [c_j;d_j]  \text{ and continuous with respect to }\lambda \notag \\ \tag{{\bf A4}}&\text{ and } \sup_{\lambda \in \Lambda}  \left\vert d\varphi _\lambda\left(x^\lambda _n\right)-  d\varphi _\lambda\left(x\right) \right\vert \xrightarrow{ \sup_{\lambda \in \Lambda}  \left\vert x_n^\lambda- x  \right\vert \rightarrow 0}0.
\end{align}
\begin{equation}
\label{hyp:tassio12_th}\tag{{\bf A5}}
 \forall 1 \leqslant j \leqslant J \quad \E \left[ \left\vert X_j \right\vert^ r \right] < \infty \text{ for some } r>4 
\end{equation}
Here $d$ is the derivation operator w.r.t. $x$, while $\partial$ will be the derivation operator w.r.t. $\lambda$.
Finally \ref{hyp:tassio12theta} becomes
\begin{equation}
\label{hyp:tassio12theta_th}\tag{{\bf A6}}
\forall  1 \leqslant j \leqslant J \quad \E \left[\sup_{\lambda \in \Lambda} \left\vert \varphi_{\lambda}\left( X_j\right)\right\vert ^ r \right]<\infty 
\text{ for some } r>4 
\end{equation}
Note that Assumption \ref{hyp:tassio12theta_th} implies that $\varepsilon$ has a moment of order $r>4$ and also that Assumption \ref{hyp:Rajput} becomes 
simpler in a parametric model which does not require a particular topology.

We impose as identifiability condition,
\begin{equation}\tag{{\bf A7}}
\label{hyp:identifTassio}
U \text{ has a unique minimizer, } \theta^\star, \text{ that belongs to the interior of }\Lambda.
\end{equation}
Note that, equivalently, this means that $\theta^*$ is the unique zero of $U$, since we are assuming that \ref{eq:H0} holds.

Now, to get sharper result about the convergence of $\inf_{\theta \in \Theta} U^n\left(\theta \right)$, one has to add the following assumptions, first on the deformation functions.
\begin{align}
\label{hyp:regul_eps}\tag{{\bf A8}}
\forall 1 \leqslant j \leqslant J &\quad \varphi^{-1}_{\theta_j^ \star} \text{ is }C^ 1 \text{ w.r.t. } x \text{ and } d\varphi^{-1}_{\theta_j^ \star} \text{ is bounded on } [a,b]\\
 &\varphi  \text{ is } C^ 2 \text{ w.r.t. } x \text{ and } \lambda \notag
\end{align}
\begin{equation}
\label{hyp:integrpartial_eps}\tag{{\bf A9}}
 \forall 1 \leqslant j \leqslant J \quad
 \E \left[ \sup_{\lambda \in \Lambda } \left\vert \partial^2\varphi_{\lambda } \left( \varphi_{\theta^\star_j}^{-1}\left(\varepsilon \right) \right) \right\vert^ 2 \right] < \infty 
\end{equation}

As said for Assumption \ref{hyp:Rajput}, the following one is more restrictive on the tail of the distribution of $\varepsilon$, excluding the Gaussian case. 
Examples of such variables with unbounded support  are given in \cite{coursTassio} p.76. Note that distributions 
with compact support and strictly positive, continuous density satisfy this assumption.
\begin{equation}\tag{{\bf A10}}\label{hyp:tassioA3}
\int_0 ^ 1 \frac{t(1-t)}{g^2\left(G^ {-1} (t) \right)}dt <\infty
\end{equation}
\subsection{Estimation of the deformation model} \label{s:s:estim}
Set $$\widehat{\theta}^n \in \arg\min_{\theta \in \Theta} U^n(\theta).$$

The results in this section are stated in the case where $\Lambda$ is a subset of $\R$. However they are still true if $\Lambda \subset \R^p$ with corresponding changes. 
The following result implies that $\widehat{\theta}^n$ is a good candidate to estimate $\theta^\star$. It is 
a simple consecuence of continuity of $U$ plus uniform convergence in probability of $U_n$ to $U$, as shown
in the proof of Theorem \ref{th:testgen}. We omit details.

\begin{prop}
\label{prop:cvestimTassio}
Under \ref{hyp:fct_th} to \ref{hyp:identifTassio}, then $$\widehat{\theta}^n\rightarrow\theta^\star\text{ in probability.}$$
\end{prop}
We can refine this result  by making the following additional assumption,
\begin{equation}
\tag{{\bf TCL}}
\label{hyp:Rj}
R_j:=\partial\varphi_{\theta^\star_j} \circ \varphi_{\theta^\star_j}^{-1} \text{ is continuous and bounded on } [a,b], \, 1\leq j\leq J.
\end{equation}
Define now $\Phi=[\Phi_{i,j}]_{1\leq i,j\leq J}$ with
\begin{equation}
\label{def:phi}
\Phi_{i,j}=-\frac{2}{J^2} \langle R_i,R_j \rangle_\mu,\, i\ne j;\quad   \Phi_{i,i}=\frac{2(J-1)}{J^2} \|R_i\|_\mu,
\end{equation}
where $\left\Vert \cdot \right\Vert_\mu$ and $\left\langle \cdot ,\cdot \right\rangle_\mu $ denote norm and inner product, respectively,
in $L^2(\mu)$. $\Phi$ is a symmetric, positive semidefinite matrix. To see this, consider $x\in\mathbb{R}^J$ and note that
\begin{eqnarray*}
x'\Phi x&=&\frac 2{J^2}\int \big(\sum_i (J-1) x_i^2 R_i^2 -2\sum_{i<j} x_i x_j R_i R_j\big) d\mu\\
&=&\frac 2{J^2}\int \sum_{i<j} (x_iR_i-x_j  R_j)^2  d\mu\geq 0.
\end{eqnarray*}
In fact, $\Phi$ is positive definite, hence invertible, unless all the $R_i$ are proportional $\mu$-a.s..
Now, we have the following Central limit Theorem, which is proved in the Appendix.
\begin{prop}
\label{prop:tcl}
Under  Assumptions \ref{hyp:fct_th} to \ref{hyp:integrpartial_eps} and \ref{hyp:Rj}, if, in addition, $\Phi$ is invertible, then
$$\sqrt{n}(\hat{\theta}^n - \theta^\star) \rightharpoonup \Phi^{-1} Y,$$
where $Y\overset d =(Y_1,\ldots,Y_J)$ with 
$$Y_{j}=\frac{2}{J} \int_0^1 R_j\circ G^{-1} \frac{\tilde{B}_j}{g \circ G^{-1} },$$
$\tilde{B_j}=B_j-\frac 1 J\sum_{k=1}^J B_k$ and $\left(B_j\right)_{1\leqslant j \leqslant J}$ independent Brownian bridges.
\end{prop}   

We note that, while, for simplicity, we have formulated Proposition \ref{prop:cvestimTassio} assuming that the deformation
model holds, a similar version can be proved (with some additional assumptions and changes in $\Phi$) in the case
when the model is false and $\theta^*$ is not the \textit{true} parameter, but the one that gives the
best (but imperfect) alignment.

\begin{remark} \label{idcond} The indentifiability condition \ref{hyp:identifTassio}
can be too strong to be realistic. Actually, for some deformation models it could
happen that $\varphi_\theta \circ \varphi_\eta=\varphi_{\theta*\eta}$ for some 
$\theta*\eta\in\Theta$. In this case, if $X_{i,j}=\varphi_{\theta_j^*}^{-1}(\varepsilon_{i,j})$ with $\varepsilon_{i,j}$ i.i.d.,
then, for any $\theta$, $X_{i,j}=\varphi_{\theta*\theta_j^*}^{-1}(\tilde{\varepsilon}_{i,j})$ with $\tilde{\varepsilon}_{i,j}=
\varphi_{\theta}({\varepsilon}_{i,j})$ which are also i.i.d. and, consequently,
$(\theta* \theta^*_1,\ldots,\theta* \theta^*_J)$ is also a zero of $U$. This applies, for instance, to location and scale models.
A simple fix to this issue is to select one of the 
signals as the reference, say the $J$-th signal, and
assume that $\theta_J^*$ is known (since it can be, in fact, chosen arbitrarily).
The criterion function becomes then $\tilde{U}(\theta_1,\ldots,\theta_{J-1})=U(\theta_1,\ldots,\theta_{J-1},\theta_J^*)$.
One could then make the (more realistic) assumption that $\tilde{\theta}^*=(\theta_1^*,\ldots,\theta_{J-1}^*)$ is the unique
zero of $\tilde{U}$ and base the analysis on $\tilde{U}_n(\theta_1,\ldots,\theta_{J-1})=U_n(\theta_1,\ldots,\theta_{J-1},\theta_J^*)$
and $\hat{\tilde{\theta}}^n=\arg\min_{\tilde{\theta}} \tilde{U}_n(\tilde{\theta})$. The results in this section can be adapted almost
verbatim to this setup. Proposition \ref{prop:tcl} holds, namely, $\sqrt{n}(\hat{\tilde{\theta}}^n-{\tilde{\theta}^*})
\rightharpoonup \tilde{\Phi}^{-1} \tilde{Y}$, with $\tilde{Y}\overset d=(Y_1,\ldots,Y_{J-1})$ and $
 \tilde{\Phi}=[\Phi_{i,j}]_{1\leq i,j\leq J-1}$. We note further that invertibility of 
$\tilde{\Phi}$ is almost granted. In fact, arguing as above, we see that
\begin{eqnarray*}
x'\tilde\Phi x=\frac 2{J^2}\int \Big( \sum_{1\leq i<j\leq J-1} (x_iR_i-x_j  R_j)^2 +  \sum_{1\leq i\leq J-1} x_i^2 R_i^2 \Big)d\mu\geq 0
\end{eqnarray*}
and $\tilde{\Phi}$ is positive definite unless $R_i=0$ $\mu$-c.s. for $i=1,\ldots,J-1$.
\end{remark}

\subsection{Asymptotic behavior of Wasserstein's variation under the null}
Here we are able to specify the speed of convergence of $\inf_{\theta \in \Theta} U^n\left(\theta \right)$ to zero
when \ref{eq:H0} holds, providing the asymptotic distribution of this statistic.

\begin{theo}
\label{th:testH0}
Under assumptions   \ref{hyp:fct_th}  to \ref{hyp:tassioA3}, \ref{hyp:Rj} and invertibility of $\Phi$,
\begin{align*}
n&\inf_{\theta \in \Theta} U_n(\theta) \rightharpoonup  
\frac 1 J\sum_{j=1}^J \int_0^1 \Big(\frac{\tilde{B}_j}{g \circ G^{-1}  }\Big)^2 -\frac 1 2 Y'\Phi^{-1}Y
\end{align*}
with $Y=(Y_1,\ldots,Y_J)$, $Y_{j}=\frac{2}{J} \int_0^1 R_j\circ G^{-1} \frac{\tilde{B}_j}{g \circ G^{-1} }$,
$\tilde{B_j}=B_j-\frac 1 J\sum_{k=1}^J B_k$ and $\left(B_j\right)_{1\leqslant j \leqslant J}$ independent Brownian bridges.
\end{theo}

A proof of Theorem \ref{th:testH0} is given in the Appendix.
As for Theorem \ref{th:testgen}, this result can be generalized to the case of different sample sizes with
straightforward changes. We also note that the result can also be adapted to the setup of Remark \ref{idcond},
replacing the correction term $\frac 1 2 Y'\Phi^{-1} Y$ by $\frac 1 2\tilde Y'\tilde\Phi^{-1} \tilde Y$.

Turning back to our goal of assessment of the deformation model \ref{eq:H0} 
based on the observed value of $\inf_{\theta \in \Theta} U^n\left(\theta \right)$,
Theorem \ref{th:testH0} gives some insight into the threshold levels for rejection
of \ref{eq:H0}. However, the limiting distribution still depends on unknown
objects and designing a tractable test requires to estimate the quantiles of this 
distribution. This will be achieved in the next section.

\section{Testing procedure with Wasserstein distance} \label{s:test}
\subsection{Bootstrap with Wasserstein distance} \label{s:boot}
In this section we present general results on Wasserstein distances that we will apply to 
estimate the asymptotic distribution of a statistic test based on an alignment with 
respect to the Wasserstein's barycenter. 
More precisely, here we consider distributions on $\R^d$ with a moment of order $r\geqslant 1$, 
that is, distributions in $\mathcal{W}_r\left(\mathbb{R}^d\right)$. ${W}_r$ will denote Wasserstein distance with $L_r$ cost, namely,
$${W}_r^r(\nu,\eta)=\inf_{\pi\in \Pi\left(\nu,\eta\right)    }\int \|y-z\|^r d\pi(y,z),$$ 
where $\left\Vert \cdot \right\Vert$ is any norm on $\R^d$. Finally, we write $\mathcal{L}(Z)$ for the law of any
random variable $Z$. 
We note 
the abuse of notation in the following, in which 
${W}_r$ is used both for Wasserstein distance on $\mathbb{R}$
and on $\mathbb{R}^d$, but this should not cause much confusion.

The next result shows that the laws of empirical transportation costs are
continuous (and even Lipschitz) functions of the underlying distributions.
\begin{prop}\label{prop:transportationcost}
Set $\nu, \nu', \eta$  probability measures in $\mathcal{W}_r\left(\mathbb{R}^d\right)$, $Y_1,\ldots,Y_n$  
i.i.d. random vectors with common law $\nu$, $Y'_1,\ldots,Y'_n$, i.i.d. with law  $\nu'$ and
write $\nu_n$, $\nu_n'$ for the corresponding empirical measures.  Then
$${W}_r(\mathcal{L}({W}_r(\nu_n,\eta)),\mathcal{L}({W}_r(\nu'_n,\eta)))\leqslant {W}_r(\nu,\nu').$$
\end{prop}
Our deformation assessment criterion concerns a particular version of the  {Wasserstein $r$-variation} of distributions 
$\nu_1,\ldots,\nu_J$ in $\mathcal{W}_r\left(\mathbb{R}^d\right)$, that we will denote in its general form by 
$$V_r(\nu_1,\ldots,\nu_J):=\inf_{\eta\in \mathcal{W}_r(\mathbb{R}^d)} \Big(\frac 1 J \sum_{j=1}^J {W}_r^r (\nu_j,\eta) \Big)^{1/r}.$$
$V_r$ is just 
the average distance to the $r$-barycenter of the set.

It is convenient to note that $V_r^r(\nu_1,\ldots,\nu_J)$
can also be expressed as
\begin{equation}
\label{eq:new_formulation}
V_r^r(\nu_1,\ldots,\nu_J)=\inf_{\pi\in \Pi(\nu_1,\ldots,\nu_J)} \int T(y_1,\ldots,y_J)d\pi(y_1,\ldots,y_J),
\end{equation}
where $\Pi(\nu_1,\ldots,\nu_J)$ denotes the set of probability measures on $\mathbb{R}^d$ with marginals
$\nu_1,\ldots,\nu_J$ and $ T(y_1,\ldots,y_J)=\min_{z\in\mathbb{R}^d} \frac 1 J \sum_{j=1}^J \|y_j-z\|^r$. A discussion about this formulation for $r=2$ and a result on existence and uniqueness of a minimizer in problem \eqref{eq:new_formulation} are given in Proposition 4.2 in \cite{agueh2010barycenters}.

Here we are interested in empirical Wasserstein $r$-variations, namely, 
the $r$-variations computed from the empirical measures
$\nu_{n_j,j}$ coming from independent samples $Y_{1, j},\ldots,Y_{n_j, j}$
of i.i.d. random variables with distribution $\nu_j$. Note that in this case
problem \eqref{eq:new_formulation} is a linear optimisation problem 
for which a minimizer always exists.

 As before, we consider 
the continuity of the law of empirical Wasserstein $r$-variations with respect to
the underlying probabilities. This is covered in the next result.

\begin{prop}\label{prop:rvariation}
With the above notation
$${W}_r^r(\mathcal{L}(V_r(\nu_{n_1,1},\ldots,\nu_{n_J,J})),\mathcal{L}(V_r(\nu'_{n_1,1},\ldots,\nu'_{n_J,J})))
\leqslant \frac 1 J \sum_{j=1}^J {W}_r^r(\nu_j,\nu'_j).$$
\end{prop}
A useful consequence of the above results is that empirical Wasserstein distances or $r$-variations can be bootstrapped
under rather general conditions. To be more precise, we take in Proposition \ref{prop:transportationcost} $\nu'=\nu_n$, the empirical
measure on $Y_1,\ldots,Y_n$ and consider a bootstrap sample $Y_1^*,\ldots,Y_{m_n}^*$ of i.i.d. (conditionally given $Y_1,\ldots,Y_n$) 
observations with common law $\nu_n$. We write $\nu_{m_n}^*$ for the empirical measure on $Y_1^*,\ldots,Y_{m_n}^*$
and $\mathcal{L}^*(Z)$ for the conditional law of $Z$ given $Y_1,\ldots,Y_n$. Proposition 
\ref{prop:transportationcost} now reads 
$${W}_r(\mathcal{L}^*({W}_r(\nu^*_{m_n},\nu)),\mathcal{L}({W}_r(\nu_{m_n},\nu)))\leqslant {W}_r(\nu_n,\nu).$$
Hence, if ${W}_r(\nu_n,\nu)=O_{\p}(1/r_n)$ for some sequence $r_n>0$ such that $r_{m_n}/r_n\to 0$ as $n\to\infty$, then,
using that ${W}_r(\mathcal{L}(aX),\mathcal{L}(aY))=a{W}_r(\mathcal{L}(X),\mathcal{L}(Y))$ for $a>0$,
we see that 
\begin{equation}
\label{eq:cv_boots1}{W}_r(\mathcal{L}^*(r_{m_n}{W}_r(\nu^*_{m_n},\nu)),\mathcal{L}(r_{m_n}{W}_r(\nu_{m_n},\nu)))\leqslant 
\frac{r_{m_n}}{r_n} r_n {W}_r(\nu_n,\nu)\to 0
\end{equation}
 in probability.\\
\indent

If in addition   $r_n{W}_r(\nu_n,\nu)\rightharpoonup \gamma \left(\nu \right)$ for a distribution $\gamma \left(\nu \right) $  then  $${r_{m_n}}{W}_r(\nu_{m_n},\nu)\rightharpoonup \gamma(\nu)$$
 which entails that if $\hat{c}_n(\alpha)$ denotes the $\alpha$ quantile of the conditional
law $\mathcal{L}^*(r_{m_n}{W}_r(\nu^*_{m_n},\nu))$  then under some regularity conditions on the distribution $\gamma(\nu)$
\begin{equation}\label{eq:bootquantile}
\p\left( r_n{W}_r(\nu_n,\nu) \leqslant \hat{c}_n(\alpha)\right)\to \alpha \quad {\rm as}\:   n\to\infty.
\end{equation}
We conclude in this case that  the quantiles of $r_n{W}_r(\nu_n,\nu)$
can be consistently estimated by the bootstrap quantiles, that is, the conditional 
quantiles of $r_{m_n}{W}_r(\nu^*_{m_n},\nu)$ (which, in turn, can be approximated through 
Monte-Carlo simulation).\\ 
As an example, if $d=1$ and $r=2$, under integrability and smoothness assumptions on $\nu$ we have 
$\sqrt{n}{W}_2(\nu_n,\nu)\rightharpoonup \left(\int_{0}^1 \frac{B^2(t)}{f^2(F^{-1}(t))} dt\right)^{1/2}$, where $f$ and $F^{-1}$ are the density and the quantile function  of $\nu$.

\subsection{Bootstrap for Wasserstein's barycenter alignement}
In the non parametric deformation model,   statistical inference  is based on the minimal Wasserstein variation
$$v^2_n:=\inf_{ \varphi \in\mathcal{G}} V^2_2(\mu^{n,1}(\varphi ),\ldots,\mu_{n,J}(\varphi ))= \inf_\mathcal{G} U_n,$$
where $\mu_{n,j}(\varphi )$ denotes the empirical measure on $Z_{1,j}(\varphi ),\ldots,Z_{n,j}(\varphi )$, where
$Z_{i,j}(\varphi )=\varphi_{j}^{-1}(X_{i,j})$ and $X_{1,j},\ldots,X_{n,j}$ are independent i.i.d. samples
from $\mu_j$. 
Consider $v_n'$, the corresponding version obtained from samples with underlying distributions $\mu_j'$, and denote 
by $\mathcal{L}\left(v_n\right)$ (reps. $\mathcal{L}\left(v'_n\right)$) the law of the random variable $v_n$ (resp. $v'_n$).

Then,  the following result holds, setting 
$\left\Vert \varphi_j\right\Vert_\infty = sup_{x \in \left(c_j;d_j\right)}\left\vert \varphi_j(x)\right\vert$. 

\begin{theo} \label{th:bootstrap}
Under  Assumption   \ref{hyp:fct}, if for all $j$ $\mathcal{G}_j \subset \mathcal{C}^1\left(c_j;d_j\right)$ and $\sup_{\varphi \in \mathcal{G}} \left\Vert  \varphi_j^\prime \right\Vert_{\infty }< \infty$, then  
$${W}_2^2(\mathcal{L}(v_n),\mathcal{L}(v'_n))
\leqslant\sup_{\varphi \in \mathcal{G}} \left\Vert  \varphi_j^\prime \right\Vert^2_{\infty } \frac 1 J \sum_{j=1}^J {W}_2^2(\mu_j,\mu'_j).$$
\end{theo}

Now consider bootstrap samples $X^*_{1,j},\ldots,X^*_{m_{n},j}$ 
of i.i.d. observations sampled from $\mu_{n,j}$, write $\mu^*_{m_n,j}$ for the empirical measure on
$X^*_{1,j},\ldots,X^*_{m_n,j}$ (conditionally to the $X_{1,j},\ldots,X_{{n},j}$)    and denote  
$V^2_2(\mu^*_{m_n,1}(\varphi ),\ldots,\mu^*_{m_n,J}(\varphi )) =  U_{m_n}^* \left( \varphi \right)$. Then we get
\begin{corollary}\label{cor:bootstrap}
If $m_n\to \infty$, and $m_n/\sqrt{n}\to 0$, then under Assumptions  \ref{hyp:fct}  to \ref{hyp:tassio12theta}, and if $\inf_\mathcal{G}U>0$,  writing $\gamma$
for the limit distribution in Theorem \ref{th:testgen}, we have that
$$\mathcal{L}^*\left( \sqrt{m_n} \left( \inf_{ \mathcal{G}} U_{m_n}^*   - \inf_{\mathcal{G}}U \right) \right) \rightharpoonup \gamma $$
in probability. In particular, if $\hat{c}_n(\alpha)$ denotes the conditional (given the $X_{i,j}$'s) $\alpha$-quantile
of $\sqrt{m_n} \left( \inf_{ \mathcal{G}} U_{m_n}^*   - \inf_{\mathcal{G}}U  \right)  $ then 
\begin{equation}
\label{eq:res_cor}
 \p\left(  \sqrt{n} \left(\inf_{\mathcal{G}} U^n -\inf_{\mathcal{G}} U \right) \leqslant \hat{c}_n(\alpha)   \right)\to \alpha. \end{equation}
\end{corollary}
Now consider the parametric deformation model and note that
the inference about it is based on the minimal Wasserstein variation
$$v^2_n:=\inf_{ \theta \in\Theta} V^2_2(\mu_{n, 1}(\theta ),\ldots,\mu_{n, J}(\theta ))= \inf_\Theta U_n,$$
where $\mu_{n, j}(\theta )$ denotes the empirical measure on $Z_{1,j}(\theta ),\ldots,Z_{n,j}(\theta )$,
$Z_{i,j}(\theta )=\varphi_{\theta_j}^{-1}(X_{i,j})$ and $X_{1,j},\ldots,X_{n,j}$ are independent i.i.d. samples
from $\mu_j$. 
We consider $v_n'$, the corresponding
version obtained from samples with underlying distributions $\mu_j'$, and denote by $\mathcal{L}\left(v_n\right)$ 
(resp. $\mathcal{L}\left(v'_n\right)$) the law of the random variable $v_n$ (resp. $v'_n$).

Then, we are able to prove the following result.

\begin{theo} \label{th:bootstrap_th}
Under  Assumptions   \ref{hyp:fct_th}, \ref{hyp:regulbis_th} and  \ref{hyp:regul_th} 
$${W}_2^2(\mathcal{L}(v_n),\mathcal{L}(v'_n))
\leqslant \sup_{x \in (c;d), \lambda \in \Lambda }\left\vert d\varphi_\lambda (x)\right \vert^2 \frac 1 J \sum_{j=1}^J {W}_2^2(\mu_j,\mu'_j).$$
\end{theo}

Now consider bootstrap samples $X^*_{1,j},\ldots,X^*_{m_{n},j}$ 
of i.i.d. observations sampled from $\mu^n_j$, write $\mu^*_{m_n,j}$ for the empirical measure on
$X^*_{1,j},\ldots,X^*_{m_n,j}$ (conditionally to the $X_{1,j},\ldots,X_{{n},j}$)    and denote  
$V^2_2(\mu^*_{m_n,1}(\theta ),\ldots,\mu^*_{m_n,J}(\theta )) =  U_{m_n}^\star \left( \theta \right)$.
\begin{corollary}\label{cor:bootstrap_th}
If $m_n\to \infty$, and $m_n/{n}\to 0$, then under Assumptions   \ref{hyp:fct_th} to \ref{hyp:tassioA3}, \ref{hyp:Rj}   and  writing $\gamma\left(G;\theta ^*\right)$
for the limit distribution in Theorem \ref{th:testH0}, we have that
$$\mathcal{L}^*\left( {m_n} \inf_{ \Theta} U_{m_n}^*    \right) \rightharpoonup \gamma\left(G; \theta^*\right)$$
in probability. In particular, if $\hat{c}_n(\alpha)$ denotes the conditional (given the $X_{i,j}$'s) $\alpha$-quantile
of $ {m_n} \inf_{\Theta} U_{m_n}^* $ then if the quantile function of $\gamma\left(G;\theta ^*\right)$ is continuous w.r.t $\alpha$ 
\begin{equation}
\label{eq:res_cor_th}
 \p\left(  {n} \inf_{\Theta} U^n \leqslant \hat{c}_n(\alpha)   \right)\to \alpha. \end{equation}
\end{corollary}

\subsection{Goodness of fit} \label{s:goodness}
In the semi parametric model, we can now provide a goodness of fit procedure. Under 
Assumptions of Theorem \ref{th:testH0} (\ref{hyp:fct_th} to \ref{hyp:tassioA3} and  \ref{hyp:Rj}) one can test the null assumption 
\begin{equation}
\label{eq:H0_test_par_v2} \inf_{\theta \in \Theta} U(\theta)=0 \tag{$\mathcal{H}_0$}
\end{equation}
versus its complementary denoted by $\mathcal{H}_1$.

In this case the test statistic is $n\inf_\Theta U_n$ and one can get the asymptotic level of a reject 
region of the form $\left\lbrace n\inf_\Theta U_n > \lambda_n \right\rbrace$ by using Corollary \ref{cor:bootstrap_th}.

More precisely,  consider bootstrap samples $X^*_{1,j},\ldots,X^*_{m_{n},j}$ 
of i.i.d. observations sampled from $\mu_{n,j}$, and write $U_{m_n}^* \left( \theta \right)$ for the corresponding criterion.
Then, if $\hat{c}_n(\alpha)$ denotes the conditional (given the $X_{i,j}$'s) $(1-\alpha)$-quantile
of $ {m_n} \inf_{\Theta} U_{m_n}^* $
\begin{equation*}
 \p\left(  {n} \inf_{\theta \in \Theta} U_n(\theta) > \hat{c}_n(\alpha)   \right)\to \alpha. \end{equation*}

Thus $\left\lbrace {n}\inf_{\theta \in \Theta} U_n(\theta)> \hat{c}_n(\alpha)   \right\rbrace$ will be a reject region 
of asymptotic level $\alpha$, and $\hat{c}_n(\alpha) $ can be computed using a Monte-Carlo method.

Note that in the case of a non parametric model, a test can be designed switching the null hypothesis. Hence set for $\Delta_0>0$ set
\begin{equation}
\label{eq:H0_test_nonpar} \inf_{\mathcal{G}}U= \Delta_0 \tag{$\mathcal{H}_0^1$}
\end{equation}
\begin{equation}
\label{eq:H1_test_nonpar} \inf_{\mathcal{G}} U<\Delta_0 \tag{$\mathcal{H}_1^1$}
\end{equation}
The test statistic in this case is $\mathcal{U}_n\left(\Delta_0\right) : =\sqrt{n} \left(\inf_{\mathcal{G}} U_n - \Delta_0 \right)$.
Then, under assumptions of  Corollary \ref{cor:bootstrap} (\ref{hyp:fct}  to \ref{hyp:tassio12theta}), 
if $\hat{c}_n(\alpha)$ denotes the conditional (given the $X_{i,j}$'s) $\alpha$-quantile
of the bootstrap version $\sqrt{m_n} \left( \inf_{ \mathcal{G}} U_{m_n}^*   - \Delta_0 \right)  $, 
under  \ref{eq:H0_test_nonpar}
\begin{equation*}
 \p\left(  \mathcal{U}_n\left(\Delta_0 \right)  \leqslant \hat{c}_n(\alpha)   \right)\to \alpha, \end{equation*}
which gives the asymptotic level of the reject region $\left\lbrace  \mathcal{U}_n\left(\Delta_0 \right) 
\leqslant \hat{c}_n(\alpha) \right\rbrace$, where $\hat{c}_n(\alpha)$ can be computed using a Monte-Carlo method.

This procedure can be made more precise under Assumptions of Theorem \ref{th:testgen} in the parametric case (\ref{hyp:fct_th} to \ref{hyp:tassio12theta_th}).  
Set for $1 \leqslant j\leqslant J$,  set $S_j \left( \theta \right)= \int_0^1  \varphi_{\theta_j}^\prime \left( F_j^{-1}\left( t \right) 
\right)\frac{B_j\left( t \right)  }{ {f_j\left(F_j^{-1}(t) \right)}} \left( \varphi_{\theta_j} \left( F_j^{-1}\left( t \right)\right) -
F^{-1}_B\left( {\theta} \right)\left( t \right) \right)dt,$ independent centered Gaussian variables. Then the result of Theorem \ref{th:testgen} can be restated as  
\begin{align*}
\sqrt{n}\left\lbrace \inf_{\Theta}  U_n -\inf_{\Theta}U  \right\rbrace \rightharpoonup \frac{2}{J}\sum_{j=1}^J S_j\left( \theta^\star \right). 
\end{align*}
Let $\sigma^2_j$ the variance of $S_j(\theta^\star)$. Set
$\widehat{L}_j\left(\frac{i}{n}\right)  = \frac{1}{2}\left\lbrace{Z_{(1)j} \left( \widehat{\theta}^n\right)}^2-{Z_{(1)j} \left( \widehat{\theta}^n\right)} ^2 \right\rbrace - \sum_{k=2}^i\left(\frac{1}{J} \sum_{p=1}^J Z_{(k)p} \left( \widehat{\theta}^n\right) \right)\left( Z_{(k)j} \left( \widehat{\theta}^n\right) -Z_{(k-1)j} \left( \widehat{\theta}^n\right)  \right)$
Then we could prove  that 
$$\widehat{\sigma}^j_n =\sum_{i= 2} ^n \frac{n-1}{n^2} \sum_{i=2}^n \widehat{L}^2_j\left(\frac{i}{n}\right)  - \frac{1}{n^2} \sum_{\substack{k ,i=2 \\ k\neq i }}^n \widehat{L} _j\left(\frac{i}{n}\right)\widehat{L} _j\left(\frac{k}{n}\right) $$
converges in probability to $\sigma^2_j$. Hence, we can now
provide a test procedure for the null assumption 
\begin{equation}
\label{eq:H0_test_par_v1} \inf_{\Theta} U \geqslant  \Delta_0 \tag{$\mathcal{H}^2_0$}
\end{equation}
versus its complementary denoted by $\mathcal{H}^2_1$.

Here we set the test statistic as $ \mathcal{V}_n\left(\Delta_0\right) : =\sqrt{n} \frac{\inf_{\Theta} U_n - \Delta_0 }{\widehat{\sigma}_n}$.
Then 
$$ \mathcal{V}_n\left(\Delta_0\right) =\sqrt{n} \frac{\inf_{\Theta} U_n -\inf_{\Theta} U }{\widehat{\sigma}_n}  +\sqrt{n} \frac{ \inf_{\Theta} U -\Delta_0 }{\widehat{\sigma}_n}$$
and if $\inf_{\Theta} U =\Delta_0 $
$$ \mathcal{V}_n\left(\Delta_0\right) \rightharpoonup Z \sim \mathcal{N}\left( 0 , 1 \right)$$
else, if $\inf_{\Theta} U >\Delta_0 $, we get that for all  $m \in \R$
$$\p\left(  \mathcal{V}_n\left(\Delta_0\right) \geqslant m \right) \rightarrow 1.$$
Then,
$$\sup_{\left( \mu_1, \dots, \mu_J \right) \text{ s.t. } \eqref{eq:H0_test_par_v1} \text{ holds}} \lim_{n\rightarrow \infty}\p\left(  \mathcal{V}_n\left(\Delta_0\right) \leqslant \lambda \right)  \leqslant \Phi\left( \lambda \right) $$
where $\Phi$ is the distribution function of the standard normal distribution. Thus we can construct a test of asymptotic level $\alpha$ by choosing the reject region $\left\lbrace \mathcal{ V}_n\left(\Delta_0\right) \leqslant  \Phi^{-1}(\alpha)\right\rbrace$.

\section{Appendix} \label{s:append} We provide here proofs of the main results in this paper. 
For those in Sections 3 and 4 our approach
relies on the consideration of quantile processes, namely,
$$\rho_{n,j}(t)=\sqrt{n}f_j(F_j^{-1}(t))(F_{n,j}^{-1}(t)-F_{j}^{-1}(t)),\quad 0<j<1, \, j=1,\ldots,J,$$
and on strong approximations of quantile processes, as in the following result that we adapt from 
\cite{csorgohorvath93} (Theorem 2.1, p. 381 there).
\begin{theo}
\label{th:csorgo}
Under \ref{hyp:tassio13}, there exist, on a rich enough probability space, inependent versions of
$\rho_{n,j}$ and independent families of Brownian bridges $\{B_{n,j}\}_{n=1}\infty$, $j=1,\ldots,J$ 
satisfying $$n^{1/2- \nu} \sup_{1/n \leqslant t \leqslant 1- 1/n} \frac{| \rho_{n,j}(t) - B_{n,j}(t) 
|}{\left( t (1-t) \right) ^\nu}= \left\lbrace \begin{array}{c} O_p (\log(n) ) 
\text{ if } \nu=0 \\ O_p (1 ) \text{ if } 0< \nu \leqslant 1/2\end{array} \right.$$
\end{theo}

We will make frequent use in this section of the following technical Lemma which 
generalizes a result in \cite{MR2435470}.
\begin{lem}
\label{lem:ln}
Under Assumption \ref{hyp:tassio12theta} 
\begin{enumerate}
\item[i)] $\sup_ {h \in \mathcal{G}_j} \sqrt{n} \int_0^\frac{1}{n}(h ( F_j^{-1}( t ) )  )^2 dt \rightarrow 0$,  
$\sup_ {h \in \mathcal{G}_j} \sqrt{n} \int_{1-\frac{1}{n}}^1(h( F_j^{-1}( t)))^2 dt \rightarrow 0$.
\item[ii)] $\sup_ {h \in \mathcal{G}_j} \sqrt{n} \int_0^\frac{1}{n}(h ( F_{n,j}^{-1}( t ) )  )^2 dt \rightarrow 0$,  
$\sup_ {h \in \mathcal{G}_j} \sqrt{n} \int_{1-\frac{1}{n}}^1(h( F_{n,j}^{-1}( t)))^2 dt \rightarrow 0$ in probability.
\item[iii)] If moreover   \ref{hyp:Rajput} holds 
\begin{equation}
\label{Acunif}
\forall k, j 
 \int_0^1 \frac{\sqrt{t(1-t)}}{f_k\left( F_k^{-1} (t) \right)} \sup_{\varphi \in \mathcal{G}}\left\vert  \varphi_j\left( F_j^{-1}\left( t \right) \right)-F^{-1}_B\left( \varphi \right)\left( t \right)  \right\vert dt < \infty
\end{equation}
\item[iv)] In the parametric case, under Assumptions \ref{hyp:regulbis_th},  \ref{hyp:tassio12theta_th}     and if $\forall k$, $F_k$ is $C^1$ with $F'_k=f_k>0$  on $(c_k,d_k)$ 
\begin{equation}
\label{Acunif_th}
\forall k, j 
 \int_0^1 \frac{\sqrt{t(1-t)}}{f_k\left( F_k^{-1} (t) \right)} \sup_{\theta \in \Theta}\left\vert  
 \varphi_{\theta_j}^{-1} \left( F_j^{-1}\left( t \right) \right)-F^{-1}_B\left( \theta \right)\left( t \right)  \right\vert dt < \infty
\end{equation}
\end{enumerate}
\end{lem}

{Our next proof is inspired  by \cite{MR2435470}. The main part concerns the study  
of $\sqrt{n} U_n(\varphi)$ uniformly in $\varphi$ in probability by using 
strong approximations of the quantile process with Brownian bridges.}

\medskip
{
\textsc{Proof of Theorem \ref{th:testgen}. }
We will work with the versions of $\rho_{n,j}$ and $B_{n,j}$ given by Theorem~\ref{th:csorgo}. We show first that
\begin{equation}
 \label{eq:step2fin}
 \sup_{\varphi  \in \mathcal{G}} \Big| \sqrt{n} \left( U_n\left(\varphi \right)  - U\left( \varphi \right) \right)- 
 \frac{1}{J} \sum_{j=1}^J S_{n,j}\left(\varphi \right) \Big|  \rightarrow 0 \text{ in probability}
 \end{equation}
with 
$S_{n,j}(\varphi) =2 \int_0^1  \varphi_j^\prime \circ F_j^{-1} (\varphi_j\circ F_j^{-1}-F_B^{-1}(\varphi)) \frac{B_{n,j}}{f_j\circ F_j^{-1}}$.
To check this we note that the fact that $\frac 1 J\sum_{j=1}^J \varphi_j\circ F_j^{-1}= F_B^{-1}(\varphi)$ and simple algebra yield 
$\sqrt{n}(U_n(\varphi)-U(\varphi))=\frac{2}{J}\sum_{j=1}^J \tilde{S}_{n,j} + \frac{1}{J} \sum_{j=1}^J\tilde{R}_{n,j} $ with 
$$
\tilde{S}_{n,j}=\sqrt{n}\int_0^1 (\varphi_j\circ F_{n,j}^{-1}-\varphi_j\circ F_j^{-1})
(\varphi_j\circ F_j^{-1}-F_B^{-1}(\varphi)),$$
$$
\tilde{R}_{n,j}=\sqrt{n}\int_0^1 [(\varphi_j\circ F_{n,j}^{-1}-\varphi_j\circ F_j^{-1}) - (F_{n,B}^{-1}(\varphi)-F_B^{-1}(\varphi))]^2.
$$
From the elementary inequality $(a_1+\cdots+a_J)^2\leq Ja_1^2+\cdots+Ja_J^2$ we get that
$$\frac{1}{J} \sum_{j=1}^J\tilde{R}_{n,j}\leq \frac{4\sqrt{n}}{J} \sum_{j=1}^J\int_0^1 (\varphi_j\circ F_{n,j}^{-1}-\varphi_j\circ F_{j}^{-1})^2$$
Now,  for every $t\in (0,1)$ we have 
\begin{equation}
\label{taylor}
\varphi_j\circ F_{n,j}^{-1}(t)- \varphi_j\circ F_{j}^{-1}(t)= 
\varphi_j^\prime ( K_{n,\varphi_j}(t))(F_{n,j}^{-1}(t)- F_{j}^{-1}(t))
\end{equation}
for some $K_{n,\varphi_j}(t)$ between $F_{n,j}^{-1}(t)$ and $F^{-1}(t)$. Assumption \ref{hyp:regul}
implies $C_j:=\sup_{\varphi_j\in\mathcal{G}_j, x\in (c_j,d_j)}|\varphi_j'(x)|<\infty$. Hence,
we have
$$\int_0^1 (\varphi_j\circ F_{n,j}^{-1}-\varphi_j\circ F_{j}^{-1})^2\leq C_j^2
\int_0^1 (F_{n,j}^{-1}-F_{j}^{-1})^2.$$
Now we can use \ref{hyp:tassio12} and argue as in the proof of Theorem 2 in \cite{MR2435470} to conclude
that $\sqrt{n}\int_0^1 (F_{n,j}^{-1}-F_{j}^{-1})^2\to 0$ in probability and, as a consequence, that
\begin{equation}
 \label{eq:intermediate}
 \sup_{\varphi  \in \mathcal{G}} \Big| \sqrt{n} \left( U_n\left(\varphi \right)  - U\left( \varphi \right) \right)- 
 \frac{1}{J} \sum_{j=1}^J \tilde{S}_{n,j}\left(\varphi \right) \Big|  \rightarrow 0 \text{ in probability}.
\end{equation}
On the other hand, the Cauchy-Schwarz's inequality shows that
\begin{eqnarray*}
\lefteqn{n\Big(\int_0^{\frac 1 n} (\varphi_j\circ F_{n,j}^{-1}-\varphi_j\circ F_j^{-1})(\varphi_j\circ F_j^{-1}-F_B^{-1}(\varphi))\Big)^2}\hspace*{1.5cm}\\
&\leq & \sqrt{n} \int_0^{\frac 1 n} (\varphi_j\circ F_{n,j}^{-1}-\varphi_j\circ F_j^{-1})^2
 \sqrt{n} \int_0^{\frac 1 n}(\varphi_j\circ F_j^{-1}-F_B^{-1}(\varphi))^2
\end{eqnarray*}
and using i) and ii) of Lemma \ref{lem:ln}, the two factors converge to zero uniformly in $\varphi$.
A similar argument works for the upper tail and allows to conclude that
we can replace in (\ref{eq:intermediate}) $\tilde{S}_{n,j}(\varphi)$ with 
$\tilde{\tilde{S}}_{n,j}(\varphi):=
2 \sqrt{n}\int_{\frac 1 n}^{1-\frac 1 n} (\varphi_j\circ F_{n,j}^{-1}-\varphi_j\circ F_j^{-1})(\varphi_j\circ F_j^{-1}-F_B^{-1}(\varphi))$.
Moreover,
\begin{eqnarray*}
\lefteqn{\sup_{\varphi  \in \mathcal{G}}  \Big| \int_0^\frac{1}{n}\varphi_j^\prime\circ F_j^{-1}\frac{B_{n,j} }{ {f_j\circ F_j^{-1}}} 
(\varphi_j \circ  F_j^{-1}-F^{-1}_B(\varphi)) \Big|}\hspace*{1cm}\\
&\leq & C_j \int_0^\frac{1}{n} \Big| \frac{B_{n,j}}
{f_j\circ F_j^{-1}}\Big| \sup_{\varphi  \in \mathcal{G}}  \big| (\varphi_j \circ  F_j^{-1}-F^{-1}_B(\varphi)) \big| 
\end{eqnarray*}
and by iii) of Lemma \ref{lem:ln} and Cauchy-Schwarz's inequality
\begin{eqnarray*}
\lefteqn{\E  \Big[ \int_0^\frac{1}{n} \Big| \frac{B_{n,j}}
{f_j\circ F_j^{-1}}\Big| \sup_{\varphi  \in \mathcal{G}}  \big| (\varphi_j \circ  F_j^{-1}-F^{-1}_B(\varphi)) \big|  \Big] } \hspace*{1cm}\\ 
& \leq & \int_0^\frac{1}{n} \frac{ \sqrt{t(1-t)}  }{ {f_j(F_j^{-1}(t))}}\sup_{\varphi  \in \mathcal{G}}  \big| \varphi_j 
( F_j^{-1}( t )) -F^{-1}_B(\varphi)
( t) \big| dt   \rightarrow 0.
\end{eqnarray*}
Hence,  $\sup_{\varphi  \in \mathcal{G}}  \Big| \int_0^\frac{1}{n}\varphi_j^\prime\circ F_j^{-1}\frac{B_{n,j} }{ {f_j\circ F_j^{-1}}} 
(\varphi_j \circ  F_j^{-1}-F^{-1}_B(\varphi)) \Big|\to 0$ in probability and similarly for the right tail. Thus (recall \eqref{taylor}), to prove 
\eqref{eq:step2fin} it suffices to show that
\begin{eqnarray}\label{finalapp}
\lefteqn{\sup_{\varphi\in\mathcal{G}}\Big|  
\int_\frac{1}{n}^{1-\frac 1 n}\varphi_j^\prime(F_j^{-1}(t))\frac{B_{n,j}(t) }{ {f_j(F_j^{-1}(t))}} 
(\varphi_j (F_j^{-1}(t))-F^{-1}_B(\varphi)(t)) dt }\hspace*{-0.4cm}\\ \nonumber
&-& \int_\frac{1}{n}^{1-\frac 1 n}\varphi_j^\prime(K_{n,\varphi_j}(t))\frac{\rho_{n,j}(t) }{ {f_j(F_j^{-1}(t))}} 
(\varphi_j (F_j^{-1}(t))-F^{-1}_B(\varphi)(t)) dt  \Big|\to 0
\end{eqnarray}
in probability. To check it we take $\nu\in(0,1/2)$ 
and use Theorem \ref{th:csorgo} to get
\begin{align}\nonumber
\int_\frac{1}{n}^{1-\frac{1}{n} } &\frac{|\rho_{n,j}(t)- B_{n,j}(t)|}{f_j(F_j^{-1}(t))}
\sup_{\varphi  \in\mathcal{G}} \big| \varphi_j (  F_j^{-1}( t )) -F^{-1}_B
( \varphi ) (t ) \big| dt \\ \label{finalapp1} &\leqslant n^{\nu -\frac{1}{2}  }O_P(1)
\int_\frac{1}{n}^{1-\frac{1}{n} } \frac{(t( 1-t)) ^\nu}{f_k(F_k^{-1}(t) )}\sup_{\varphi  \in\mathcal{G}} \big| 
\varphi_j (  F_j^{-1}( t )) -F^{-1}_B( \varphi )( t) \big|  dt\to 0
\end{align}
in probability (using dominated convergence and iii) of Lemma \ref{lem:ln}).

We observe next that, for all $t\in (0,1)$, 
$\sup_{ \varphi_j \in  \mathcal{G}_j}| K_{n,\varphi_j}(t) - F_j^{-1}(t)| \rightarrow 0$ almost surely, since
$ K_{n,\varphi_j}(t)$ lies between $F_{n,j}^{-1}(t)$ and $F_{j}^{-1}(t)$. Therefore, using Assumption \ref{hyp:regul} 
we see that 
$\sup_{ \varphi_j\in  \mathcal{G}_j }| \varphi_j^\prime ( K_{n,\varphi_j}(t) ) - \varphi_j^\prime (F_j^{-1}(t)| \rightarrow 0$ almost surely
while, on the other hand $\sup_{ \varphi_j\in  \mathcal{G}_j }| \varphi_j^\prime ( K_{n,\varphi_j}(t) ) - \varphi_j^\prime (F_j^{-1}(t))|\leq 2C_j$.
But then, by dominated convergence we get that
$$ \E \Big[\sup_{ \varphi_j\in  \mathcal{G}_j }| \varphi_j^\prime ( K_{n,\varphi_j}(t) ) - \varphi_j^\prime (F_j^{-1}(t)) |^2 \Big]  \rightarrow 0.
$$
Since by iii) of Lemma \ref{lem:ln} we have that $t\mapsto \frac {\sqrt{t(1-t)}}{f_j(F_j^{-1}(t))} \sup_{\varphi\in\mathcal{G}}|\varphi_j(F_j^{-1}(t))-
F_B^{-1}(\varphi)(t) |$ is integrable we conclude that 
$$\E \sup_{\varphi\in \mathcal{G}} \int_{\frac1 n}^{1-\frac 1 n}| \varphi_j^\prime ( K_{n,\varphi_j}(t) ) - \varphi_j^\prime (F_j^{-1}(t)) |
\frac {|B_{n,j}(t)|}{f_j(F_j^{-1}(t))} |\varphi_j(F_j^{-1}(t))-
F_B^{-1}(\varphi)(t) |dt$$
tends to 0 as $n\to\infty$ and, consequently,
$$\sup_{\varphi\in \mathcal{G}} \int_{\frac1 n}^{1-\frac 1 n}| \varphi_j^\prime ( K_{n,\varphi_j}(t) ) - \varphi_j^\prime (F_j^{-1}(t)) |
\frac {|B_{n,j}(t)|}{f_j(F_j^{-1}(t))} |\varphi_j(F_j^{-1}(t))-
F_B^{-1}(\varphi)(t) |dt$$
vanishes in probability. Combining this fact with (\ref{finalapp1}) we prove (\ref{finalapp}) and, as a consequence, (\ref{eq:step2fin}).

}

\medskip
Observe now that for all $n \in \N$, $(S_{n,j}(\varphi))_{1 \leqslant j \leqslant J }$ has the same law as $\left(  S_j\left( \varphi \right)\right)_{1 \leqslant j \leqslant J }$
with 
$$S_{j}(\varphi) =2 \int_0^1  \varphi_j^\prime \circ F_j^{-1} (\varphi_j\circ F_j^{-1}-F_B^{-1}(\varphi)) \frac{B_{j}}{f_j\circ F_j^{-1}}$$
and $\left(B_j\right)_{1 \leqslant j \leqslant J}$ independent standard Brownian bridges.
Set $S = \frac{1}{J}\sum_{j=1}^J S_j$.
Now, \eqref{eq:step2fin} implies that 
\begin{equation}
\label{eq:cvloiS}\sqrt{n} \left( U^n\left( \cdot \right)  - U\left( \cdot \right) \right) \rightharpoonup S\left(\cdot \right)
\end{equation}
in the space $  L^\infty \left(\mathcal{G} \right)$ (we denote by $\left\Vert \cdot \right\Vert_\infty$ the norm on this space).
%
%
%
%
%
%
%
%
%
From Skohorod Theorem we know that there exists some probability space on which the convergence 
\eqref{eq:cvloiS} holds almost surely. From now on, we place us on this space.
Then, for  $\varphi, \rho \in \mathcal{G}$
\begin{eqnarray*}
\lefteqn{|S_j(\varphi)-S_j(\rho)|
\leq 2 \sup_{(c_j,d_j)}|\varphi_j'-\rho_j'| \Big|\int_0^1 \frac{B_{j}}{f_j\circ F_j^{-1}} (\varphi_j\circ F_j^{-1}-F_B^{-1}(\varphi)) \Big|
} \hspace*{1cm}\\ 
&&+2 \Big|\int_0^1 \frac{B_{j}}{f_j\circ F_j^{-1}}\rho_j'\circ F_j^{-1} (\varphi_j\circ F_j^{-1}-\rho_j\circ F_j^{-1}) \Big|\\
&&\leq 2 \sup_{(c_j,d_j)}|\varphi_j'-\rho_j'| \sup_{\varphi\in\mathcal{G}}\Big|\int_0^1 \frac{B_{j}}{f_j\circ F_j^{-1}} (\varphi_j\circ F_j^{-1}-F_B^{-1}(\varphi)) \Big|\\
&&+2\sup_{(c_j,d_j)} |\rho_j^\prime |  \Big({\textstyle\int_0^1 \big| \frac{B_{j}}{f_j\circ F_j^{-1}}}\big|^{q} \Big)^{1/q}
\Big({\textstyle  \int_0^1 |\varphi_j\circ F_j^{-1}-\rho_j\circ F_j^{-1}|^{p_0}
}\Big)^{1/p_0}
\end{eqnarray*}
But using iii) of Lemma \ref{lem:ln} 
\begin{eqnarray*}
\lefteqn{\E  \Big[ \sup_{\varphi\in\mathcal{G}}\Big|\int_0^1 \frac{B_{j}}{f_j\circ F_j^{-1}} (\varphi_j\circ F_j^{-1}-F_B^{-1}(\varphi))\Big|\Big]}\hspace*{0.5cm} \\ 
&&\leqslant \int_0^1\frac{ \sqrt{t(1-t)}  }{ {f_k\left(F_k^{-1}(t) \right)}}\sup_{\varphi  \in \mathcal{G}}  \left\vert \varphi_j \left( F_j^{-1}\left( t \right)\right) -F^{-1}_B\left( \varphi \right)\left( t \right) \right\vert dt < \infty
\end{eqnarray*}
Hence, almost surely, $\sup_{\varphi\in\mathcal{G}}\Big|\int_0^1 \frac{B_{j}}{f_j\circ F_j^{-1}} (\varphi_j\circ F_j^{-1}-F_B^{-1}(\varphi))\Big|<\infty$.
Furthermore, from Assumption \ref{hyp:Rajput}, we get that a.s. 
$$\int_0^1 \Big({\textstyle \frac{B_{j}}{f_j\circ F_j^{-1}}}\Big)^{q}< \infty $$
and thus, for some random variable $T$ a.s. finite , and  $\varphi$, $\rho \in \mathcal{G}$, we get 
$$\left\vert S_j \left( \varphi \right) - S_j \left( \rho\right) \right\vert
\leqslant T \left\Vert \varphi - \rho \right\Vert_{\mathcal{G}}.$$
Thus, we deduce that  $\left(S_j\right)_{1\leqslant j \leqslant J}$ are almost surely continuous 
functions on  $\mathcal{G}$, endowed with the norm  $\left\Vert \cdot \right\Vert_{\mathcal{G}}$.

Observe now that
\begin{equation}
\label{ineq:bornesup}
\sqrt{n} \Big(\inf_{\mathcal{G}} U_n - \inf_{\mathcal{G}} U \Big)\leq \sqrt{n} 
\inf_\Gamma U_n - \sqrt{n} \inf_\Gamma  U = \inf_\Gamma \sqrt{n} \left( U_n - U \right).
\end{equation}
On the other hand, if we consider the  (a.s.) compact set $\Gamma_n = \lbrace  \varphi \in {\mathcal{G}}:\, U\left( \varphi  \right) \leqslant 
\inf_{{\mathcal{G}}} U + 
\frac{2}{\sqrt{n}} \left\Vert \sqrt{n} \left( U_n - U \right) \right\Vert_\infty \rbrace$, then, if $\varphi \notin \Gamma^n$,
$$U_n\left(\varphi  \right)  \geqslant \inf_{\mathcal{G}} U + {2} \left\Vert  \left( U_n - U \right) \right\Vert_\infty - \left\Vert  \left( U_n - U \right) 
\right\Vert_\infty, $$
which implies 
$$U_n\left( \varphi \right)\geqslant \inf_{\mathcal{G}} U + \left\Vert  \left( U_n - U \right) \right\Vert_\infty, $$
while if $\varphi \in \Gamma$, then,
$$U_n\left( \varphi \right) = \inf_{\mathcal{G}} U + U^n\left( \varphi \right) - 
U \left( \varphi \right) \leqslant \inf_{\mathcal{G}} U + \left\Vert  \left( U_n - U \right) \right\Vert_\infty.$$
Thus, necessarily, $ \inf_ { \mathcal{G}} U_n=\inf_ { \Gamma_n} U_n=\inf_ { \Gamma_n} (U_n-U+U)\geq \inf_ { \Gamma_n} (U_n-U)+\inf_{\Gamma_n} U=
\inf_ { \Gamma_n} (U_n-U)+\inf_{\Gamma} U$. Together with \eqref{ineq:bornesup} this entails
\begin{equation}
\label{eq:encadrement}
 \inf_{\Gamma_n} \sqrt{n} ( U_n -  U )\leqslant \sqrt{n} \big(\inf_{\mathcal{G}} U_n - \inf_{\mathcal{G}} U\big)  \leqslant \inf_{\Gamma} \sqrt{n} \left( U_n -  U \right).
\end{equation}
Note that for the versions that we are considering $\|\sqrt{n}(U_n-U)-S\|_\infty\to 0$ a.s.. In particular, this implies
that $\inf_{\Gamma} \sqrt{n} \left( U_n -  U \right)\to \inf_{\Gamma} S$ a.s.. Hence, the proof will be complete if we 
show that a.s.
\begin{equation}\label{finalpart}
\inf_{\Gamma_n} \sqrt{n} \left( U_n -  U \right)\to \inf_{\Gamma} S.
\end{equation}
To check this last point, consider a sequence $\varphi_n\in\Gamma_n$ such that 
$ \sqrt{n} ( U_n(\varphi_n) -  U(\varphi_n))\leq \inf_{\Gamma_n} \sqrt{n} (U_n -  U )+\frac 1 n$.
By compactness of $\mathcal{G}$, taking subsequences if necessary, $\varphi_n\to \varphi_0$ for some $\mathcal{G}$.
Continuity of $U$ yields $U(\varphi_n)\to U(\varphi_0)$ and as a consequence, that $U(\varphi_0)\leq \inf_{\mathcal{G}} U$,
that is, $\varphi_0\in\Gamma$ a.s.. 
Furthermore,   
\begin{eqnarray*}
\lefteqn{\big| \sqrt{n} (U_{n}-U)(\varphi_n) - S(\varphi_0) \big|}\hspace*{2cm}\\
&&
\leqslant \left\Vert  \sqrt{n} \left(U_{n}-U \right)- S  \right\Vert_\infty + 
\left\vert S\left(\varphi_n\right)- S \left( \varphi_0 \right)\right\vert \to 0.
\end{eqnarray*}
This shows that
\begin{equation} \label{eq:b_sup}
\liminf \inf_{\Gamma_n}\sqrt{n} \left(U_n-U \right) \geq S\left( \varphi_0 \right)  \geqslant  \inf_\Gamma S
\end{equation}
and yields (\ref{finalpart}). This completes the proof.

\hfill $\Box$

\medskip

\textsc{Proof of Proposition \ref{prop:tcl}.}
We denote by
$\partial_j$ the derivative operator w.r.t. $\theta_j$, $1\leqslant j \leqslant J$
and $\partial_{j,k}$ for second order partial derivatives. We note that $\mathcal{H}$
entails that the empirical
d.f. on the $j$-th sample, $F_{n,j}(t)$, satisfies $F_{n,j}(t)=G_{n,j}(\varphi_{\theta_j^*}(t))$
with $G_{n,j}$ the empirical d.f. on the $\varepsilon_{i,j}$'s (which are i.i.d. $\mu$, with d.f. $G$).
We write now $\rho_{n,j}$ for the quantile process based on the $\varepsilon_{i,j}$'s. We write $B_{n,j}$
for independent Brownian bridges as given by
Theorem \ref{th:csorgo} (observe that \eqref{hyp:tassio13_eps} grants the existence of such $B_{n,j}$'s).

Assumption  \ref{hyp:Rj} implies that  $\partial\varphi_{\theta^\star_j}\in  L^2 (X_j)$. 
Moreover,  with Assumptions \ref{hyp:regul_eps}, \ref{hyp:integrpartial_eps} and compactness of 
$\Theta$, we deduce that $\sup_{\lambda \in \Lambda }\partial\varphi_{\lambda}\in  L^2 (X_j)$.
On the other hand, since $\varepsilon$ has a moment of order $r>4$,  arguing
as in the proof of point 3 in Lemma \ref{lem:ln} we have that 
\begin{equation}
\label{lem3ln}
\int_0^1 \frac{\sqrt{t(1-t)}}{g \left( G^{-1} (t) \right)}dt < \infty.
\end{equation}
From \ref{hyp:regul_eps} and \ref{hyp:integrpartial_eps} we have that $U_n$ is a $C^2$ function
and derivatives can be omputed by differentiation under the integral sign. This implies that
\begin{eqnarray}\nonumber
\partial_j U_n \left( \theta \right) = \frac{2}{J} \int_0^1 \partial\varphi_{\theta_j} (F_{n,j}^{-1}(t) )
\Big(\varphi_{\theta_j} (F_{n,j}^{-1}(t) ) -{ \frac 1 J}\sum_{k=1}^J \varphi_{\theta_k} (F_{n,k}^{-1}(t))\Big) dt,\\
\label{eq:derivatives}
\partial^2_{p,q}  U_n ( \theta ) =- \frac{2}{J^2} \int_0^1 \partial\varphi_{\theta_p}
( F_{n,p}^{-1}( t ))\partial\varphi_{\theta_q} (F_{n,q}^{-1}( t ) )  dt, \quad p\ne q
\end{eqnarray}
and
\begin{eqnarray*}
\partial^2_{p,p}  U_n (\theta) &=&\frac{2}{J} \int_0^1 \partial^2\varphi_{\theta_p} ((F_{n,p})^{-1}(t))
\Big(\varphi_{\theta_p} (F_{n,p}^{-1}(t) ) -{ \frac 1 J}\sum_{k=1}^J \varphi_{\theta_k} (F_{n,k}^{-1}(t))\Big)\\
&&+
\frac{2( J-1) }{J^2} \int_0^1 (\partial\varphi_{\theta_p}  (F_{n,p}^{-1}(t))^2 dt.  
\end{eqnarray*}
Similar expressions are obtained for the derivatives of $U(\theta)$ (replacing everywhere
$F_{n,j}^{-1}$ with $F_{j}^{-1}=\varphi_{\theta_j^*}^{-1}\circ G^{-1}$). We write
$DU_n(\theta)=(\partial_j U_n(\theta))_{1\leq j\leq J}$, $DU(\theta)=(\partial_j U(\theta))_{1\leq j\leq J}$
for the gradients and $\Phi_n(\theta)=[\partial^2_{p,q}U_n ( \theta )]_{1\leq p,q\leq J}$,
$\Phi(\theta)=[\partial^2_{p,q}U (\theta )]_{1\leq p,q\leq J}$for the Hessians of $U_n$
and $U$. Note that $\Phi^*=\Phi(\theta^*)$ is assumed to be invertible.

Recalling that $R_j= \partial \varphi_{\theta_j^\star} \circ \varphi_{\theta_j^ \star}^{-1}$, from the fact
 $DU(\theta^*)=0$ we see that
\begin{equation}
\sqrt{n} \partial_j U_n( \theta^\star)  =  \frac{2}{J} 
\int_0^1   R_j(G_{n,j}^{-1}(t))\frac{\rho_{n,j}(t)-\frac 1 J\sum_{k=1}^J\rho_{n,k}(t) }{g ( G^{-1} (t)) }dt.\label{eq:expr_D}
\end{equation}
Now, using Assumption \ref{hyp:Rj} and arguing as in the proof of Theorem \ref{th:testgen} we conclude
that 
\begin{align*}
\Bigg|\int_0^{1}  R_j (G_{n,j}^{-1}(t))\frac{\rho_{n,k}(t)}{g \left( G^{-1} (t) \right) }
dt -\int_0^{1} R_j \left( G^{-1}\left( t \right)\right) \frac{B_{n,k}(t)}{g \left( G^{-1} (t) \right) } dt\Bigg|  \rightarrow 0
\end{align*}
in probability and, consequently, 
\begin{equation}\label{gradconvergence}
\Bigg|\sqrt{n} \partial_j U_n( \theta^\star)
-\frac 2 J\int_0^{1} R_j \left( G^{-1}\left( t \right)\right) \frac{B_{n,j}(t)-\frac 1 J\sum_{k=1}^J B_{n,k}(t)}{g \left( G^{-1} (t) \right) } dt\Bigg|  \rightarrow 0
\end{equation}
in probability.

A Taylor expansion of $\partial_j U_n$ around $\theta^*$ shows that for some 
$\tilde{\theta }^n_{j}$ between  $\hat{\theta}^n$ and $\theta^\star$ we have
$$\partial_jU_n( \hat{\theta }^n)= \partial_jU^n( \theta ^\star) +
(\partial_{1j}^2U_n (\tilde{\theta }^n_{j}), \ldots, \partial_{Jj}^2U_n 
(\tilde{\theta }^n_{j}) )\cdot( \hat{\theta}^n - \theta^\star )$$
and because $\hat{\theta}^n$ is a zero of $DU_n$, we obtain
$$-\partial_jU^n( \theta ^\star)=
(\partial_{1j}^2U_n (\tilde{\theta }^n_{j}), \ldots, \partial_{Jj}^2U_n 
(\tilde{\theta }^n_{j}) )\cdot( \hat{\theta}^n - \theta^\star ).$$
Writing $\tilde{\Phi}_n$ for the $(J-1)\times (J-1)$ matrix whose $J-1$-th
row equals $(\partial_{1j}^2U_n (\tilde{\theta }^n_{j}), \ldots, \partial_{Jj}^2U_n 
(\tilde{\theta }^n_{j}) )$, $j=2,\ldots,J$, we can rewrite the last expansion as
\begin{equation}
\label{eq:eqbase}
- \sqrt{n}DU_n( \theta^*) =\tilde{\Phi}_n \sqrt{n}(\hat{\theta}^n - \theta^\star).
\end{equation}

We show next that $\tilde{\Phi}_n\to \Phi^*=\Phi(\theta^*)$ in probability. Recalling \eqref{eq:derivatives},
we consider first $\int_0^1 (\partial\varphi_{\tilde{\theta }^n_p}(F_{n,p}^{-1}( t)))^2 dt$. We have
\begin{eqnarray*}
\lefteqn{\Big(\int_0^1 (\partial\varphi_{\tilde{\theta }^n_p}(F_{n,p}^{-1}( t))-\partial\varphi_{{\theta }^*_p}(F_{p}^{-1}( t)))^2 dt 
\Big)^{1/2}}\hspace*{1cm}\\
&\leq& 
\Big(\int_0^1 (\partial\varphi_{\tilde{\theta }^n_p}(F_{n,p}^{-1}( t))-\partial\varphi_{{\theta }^*_p}(F_{n,p}^{-1}( t)))^2 dt\Big)^{1/2}\\
&&+
\Big(\int_0^1 (\partial\varphi_{{\theta }^*_p}(F_{n,p}^{-1}( t))-\partial\varphi_{{\theta }^*_p}(F_{p}^{-1}( t)))^2 dt 
\Big)^{1/2}\\
&\leq &
\Big(\int_0^1 \sup_{\lambda\in\Lambda}\big|\partial ^2 \varphi_{\lambda}(F_{n,p}^{-1}( t))\big|^2dt\Big)^{1/2}  |\tilde{\theta}_p^n-\theta_p^*|\\
&&+
\Big(\int_0^1 (R_p(G_{n,p}^{-1}(t))-R_p(G_{p}^{-1}(t)))^2dt\Big)^{1/2}\to 0
\end{eqnarray*}
in probability, where we have used assumptions \ref{hyp:integrpartial_eps}, \ref{hyp:Rj} and Proposition \ref{prop:cvestimTassio}.
A similar argument shows that $\int_0^1 (\varphi_{\tilde{\theta }^n_p}(F_{n,p}^{-1}( t))-\varphi_{{\theta }^*_p}(F_{p}^{-1}( t)))^2 dt$
in probability.
As a consequence, we conclude
\begin{equation}\label{eq:eqbase2}
\tilde{\Phi}_n\to \Phi^*,\quad \mbox{in probability.}
\end{equation}
Now, \eqref{eq:eqbase}, \eqref{gradconvergence} \eqref{eq:eqbase2} together with Slutsky's 
Theorem complete the proof.

\hfill $\Box$

\medskip
\textsc{Proof of Theorem \ref{th:testH0}.} We consider the same notation and setup as in the proof of 
Proposition \ref{prop:tcl}. Since $DU_n(\hat{\theta}^n)=0$, a Taylor expansion around
$\hat{\theta}^n$ shows that 
\begin{equation}\label{secondorderTaylor}
nU_n(\theta^*)-nU_n(\hat{\theta}^n)=\frac 1 2 (\sqrt{n} (\hat{\theta}^n-\theta^*))' \Phi(\tilde{\theta}_n)
(\sqrt{n} (\hat{\theta}^n-\theta^*))
\end{equation}
for some $\tilde{\theta}_n$ between $\hat{\theta}^n$ and $\theta^*$. Arguing as in the proof of Proposition
\ref{prop:tcl} we see that $\Phi(\tilde{\theta}_n)\to \Phi^*$ in probability. Hence, to complete the
proof if suffices to show that 
$$n U_n(\theta^*)-\frac 1 J \sum_{j=1}^k \int_0^1 \frac {\big(B_{n,j}(t)-\frac 1 J\sum_{k=1}^JB_{n,k}(t) \big)^2}{g(G^{-1}(t))^2}dt \to 0$$
in probability. Since
$$n U_n(\theta^*)=\frac 1 J \sum_{j=1}^k \int_0^1 \frac {\big(\rho_{n,j}(t)-\frac 1 J\sum_{k=1}^J\rho_{n,k}(t) \big)^2}{g(G^{-1}(t))^2}dt,$$
this amounts to proving that 
$$ \int_0^1 \frac {\big(\rho_{n,j}(t)-B_{n,j}(t) \big)^2}{g(G^{-1}(t))^2}dt\to 0$$
in probability.

Taking $\nu\in(0,\frac 1 2)$ in Theorem \ref{th:csorgo} we see that 
$$ \int_{\frac 1 n}^{1-\frac 1 n} \frac {\big(\rho_{n,j}(t)-B_{n,j}(t)) \big)^2}{g(G^{-1}(t))^2}dt\leq O_P(1) \frac 1{n^{1-2\nu}}
\int_{\frac 1 n}^{1-\frac 1 n} \frac {(t(1-t))^{2\nu}}{g(G^{-1}(t))^2}\to 0,$$
using condition \eqref{hyp:tassioA3} and dominated convergence. From \eqref{hyp:tassioA3}
we also see that $ \int_{1-\frac 1 n}^1 \frac {B_{n,j}(t)^2}{g(G^{-1}(t))^2}dt\to 0$ in probability. 
Condition \eqref{hyp:tassioA3} implies also that $ \int_{1-\frac 1 n}^1 \frac {\rho_{n,j}(t)^2}{g(G^{-1}(t))^2}dt\to 0$
in probability,
see \cite{SamworthJohnson}. Similar considerations apply to the left tail and complete the proof.

\hfill $\Box$

\medskip
\textsc{Proof of Proposition \ref{prop:transportationcost}.}
We set $T_n
={W}_r(\nu_n, \eta)$ and 
$T'_n={W}_r(\nu'_n, \eta)$ and $\Pi_n(\eta)$ for the set of probabilities on $\{1,\ldots,n\}\times \mathbb{R}^d$
with first marginal equal to the discrete uniform distribution on $\{1,\ldots,n\}$ and second marginal equal to
$\eta$ and note that we have $T_n=\inf_{\pi\in\Pi_n(\eta)} a(\pi)$ if we denote
$$a(\pi)= \left(\int_{\{1,\ldots,n\}\times \mathbb{R}^d} \|Y_i-z\|^rd\pi(i,z)\right)^{1/r}.$$
We define similarly $a'(\pi)$ from the $Y'_i$ sample to get $T'_n=\inf_{\pi\in\Pi_n(\eta)} a'(\pi)$. But
then, using the inequality $ \left\vert \left\Vert a \right\Vert -\left\Vert b \right\Vert \right\vert \leqslant \left\Vert a-b \right\Vert $, 
\begin{equation*}
|a(\pi)-a'(\pi)|\leqslant  \left(\int_{\{1,\ldots,n\}\times \mathbb{R}^d} \|Y_i-Y_i'\|^rd\pi(i,z)\right)^{1/r}=
\left(\frac 1 n\sum_{i=1}^n\|Y_i-Y_i'\|^r\right)^{1/r} \label{eq:boots1}
\end{equation*}
This implies that
\begin{equation*}|T_n-T_n'|^r\leqslant \frac 1 n\sum_{i=1}^n\|Y_i-Y_i'\|^r.\label{eq:boots2}
\end{equation*}

If we take now $(Y,Y')$ to be an optimal coupling of $\nu$ and $\nu'$, so that $E\left[ \|Y-Y'\|^r \right]={W}_r^r(\nu,\nu')$
and $(Y_1,Y'_1),\ldots,(Y_n,Y'_n)$ to be i.i.d. copies of $(Y,Y')$ we see that
for the corresponding realizations of $T_n$ and $T_n'$ we have
$$\E\left[|T_n-T_n'|^r\right]\leqslant \frac 1 n\sum_{i=1}^n \E\left[\|Y_i-Y_i'\|^r\right]={W}_r(\nu,\nu')^r.$$
But this shows that ${W}_r(\mathcal{L}(T_n),\mathcal{L}(T_n'))\leqslant {W}_r(\nu,\nu')$, as claimed.

\hfill $\Box$

\textsc{Proof of Proposition \ref{prop:rvariation}.}
We write $V_{r,n}=V_r(\nu_{n_1,1},\ldots,\nu_{n_J,J})$ and $V'_{r,n}=V_r(\nu'_{n_1,1},\ldots,\nu'_{n_J,J})$.
We note that
$$V_{r,n}^r=\inf_{\pi\in \Pi(U_1,\ldots,U_J)} \int T(i_1,\ldots,i_J) d\pi(i_1,\ldots,i_J),$$
where $U_j$ is the discrete uniform distribution on $\{1,\ldots,n_j\}$ and $T(i_1,\ldots,i_J)=\min_{z\in \mathbb{R}^d} 
\frac 1 J \sum_{j=1}^J\|Y_{i_j,j}$ $-z\|^r$. We write $T'(i_1,\ldots,i_J)$ for the equivalent function computed from the
$Y'_{i,j}$'s. Hence we have
$$|T(i_1,\ldots,i_J)^{1/r}-T'(i_1,\ldots,i_J)^{1/r}|^r\leqslant\frac 1 J \sum_{j=1}^J \|Y_{i_j,j}-Y'_{i_j,j}\|^r,$$
which implies
\begin{eqnarray*}
\lefteqn{\left|\left(\int T(i_1,\ldots,i_J) d\pi(i_1,\ldots,i_J)  \right)^{1/r}- \left(\int T(i_1,\ldots,i_J) 
d\pi(i_1,\ldots,i_J)  \right)^{1/r}\right|^r}\hspace*{0.5cm} \\
&\leqslant&  \int\frac 1 J \sum_{j=1}^J \|Y_{i_j,j}-Y'_{i_j,j}\|^r d\pi(i_1,\ldots,i_J)\\
&=&\frac 1 J \sum_{j=1}^J\int \| Y_{i_j,j}-Y'_{i_j,j}\|^rd\pi(i_1,\ldots,i_J)
=\frac 1 J \sum_{j=1}^J\left( \frac 1 {n_j} \sum_{i=1}^{n_j} \| Y_{i,j}-Y'_{i,j}\|^r\right)
\end{eqnarray*}
So,
$$|V_{r,n}-V'_{r,n}|^r\leqslant 
\frac 1 J \sum_{j=1}^J\left( \frac 1 {n_j} \sum_{i=1}^{n_j} \| Y_{i,j}-Y'_{i,j}\|^r\right).$$
If we take $(Y_j,Y'_j)$ to be an optimal coupling of $\nu_j$ and $\nu_j'$ and $(Y_{1,j},Y'_{1,j}),\ldots,$
$(Y_{n_j,j},Y'_{n_j,j})$ to be i.i.d. copies of $(Y_j,Y'_j)$, for $j=1,\ldots,J$, then
we obtain 
$$\E\left[|V_{r,n}-V'_{r,n}|^r \right]\leqslant 
\frac 1 J \sum_{j=1}^J\left( \frac 1 {n_j} \sum_{i=1}^{n_j} \E\left[\| Y_{i,j}-Y'_{i,j}\|^r\right]\right)=\frac 1 J \sum_{j=1}^J {W}_{r}^r(\nu_j,\nu'_j).$$
The conclusion follows.

\hfill $\Box$

\medskip
\textsc{Proof of Theorem \ref{th:bootstrap}.}
We can mimic the argument in the proof of Proposition \ref{prop:rvariation} to get
an upper bound on the Wasserstein distance between the laws of $v_n$ and $v_n'$, the corresponding
version obtained from samples with underlying distributions $\mu_j'$. In fact, arguing as above, we can write
$$v_n^2=\inf_{\varphi \in\mathcal{G}} \left[ \inf_{ \pi\in \Pi(U_1,\ldots,U_J)}  
\int T(\varphi ;i_1,\ldots,i_J) d\pi(i_1,\ldots,i_J)\right],$$
where $T(\varphi  ;i_1,\ldots,i_J)=\min_{y\in\mathbb{R}} \frac 1 J \sum_{j=1}^J (Z_{i_j,j}(\varphi )-y)^2$.
We write $T'(\varphi  ;i_1,$ $\ldots,i_J)$ for the same function computed on the $Z'_{i,j}(\varphi  )$'s 
and set $$\|\varphi^\prime \|_\infty:=\sup_{\substack{ x \in (c;d)  \\ \varphi \in \mathcal{G} }} |\varphi^\prime _j(x)|.$$
 Now, from the  fact $(Z_{i,j}(\varphi)-Z'_{i,j}(\varphi))^2\leqslant \|\varphi^\prime \|_\infty^2 (X_{i,j}-X'_{i,j})^2$ 
we see that
$$|T(\varphi ;i_1,\ldots,i_J)^{1/2}-T'(\varphi ;i_1,\ldots,i_J)^{1/2}|^2\leqslant 
\|\varphi^\prime  \|_\infty^2\frac 1 J \sum_{j=1}^J (X_{i_j,j}-X'_{i_j,j})^2$$
and, as a consequence, that 
\begin{eqnarray*}
\lefteqn{| V_2\left(\mu _{ 1}^n(\varphi ),\ldots,\mu ^n_{J }(\varphi )\right) -V_2\left(\mu^{\prime \ n}_{1}
(\varphi ),\ldots,\mu^{\prime \ n} _{J}(\varphi )\right)| ^2}\hspace*{3cm}\\
&\leqslant &\frac 1 J \sum_{ j =1}^J \sum_ {i_j=1}^{n_j} \frac{1}{n_j} 
\|\varphi^\prime  \|_\infty^2(X_{i_j,j}-X'_{i_j,j})^2
\end{eqnarray*}
and then 
$$(v_n-v'_n)^2\leqslant \|\varphi^\prime  \|_\infty^2 \frac 1 J \sum_{j=1}^J \left({\textstyle \frac 1 {n_j}\sum_{i=1}^{n_j}(X_{i,j}-X'_{i,j})^2}
\right).$$
If, as in the proof of Proposition \ref{prop:rvariation}, we assume that $(X_{i,j},X'_{i,j})$, $i=1,\ldots,n_j$ are i.i.d. copies of an
optimal coupling for $\mu_j$ and $\mu_j'$, with different samples independent from each other we obtain that
$$\E\left[ (v_n-v'_n)^2 \right]\leqslant \|\varphi^\prime  \|_\infty^2 \frac 1 J \sum_{j=1}^J {W}_2^2(\mu_j,\mu_j').$$
\hfill $\Box$

\medskip
\textsc{Proof of Corollary \ref{cor:bootstrap}.}
In Theorem \ref{th:bootstrap}, take $\mu'_j = \mu_{n,j} $, and set 
$v_{m_n}^{*}:=\inf_{ \varphi \in\mathcal{G}} 
V_2(\mu^*_{m_n,1}(\varphi ),\ldots,\mu^*_{m_n,J}(\varphi ))$. Then, conditionally to the $X_{1,j},\ldots,X_{{n},j}$, the result of Theorem \ref{th:bootstrap} reads now
 \begin{equation*} {W}_2^2(\mathcal{L}(v_{m_n}),\mathcal{L}(v_{m_n}^{*} ))
\leqslant \sup_{\varphi \in \mathcal{G} }\left\Vert \varphi^\prime_j\right \Vert_{\infty}^2 \frac 1 J \sum_{j=1}^J {W}_2^2(\mu_j,\mu_{n,j}).\end{equation*}
Now, let $v^2:=\inf_{ \varphi \in\mathcal{G}} 
M \left( \varphi \right)$. Then, 
\begin{eqnarray} \label{eq:cor1}{W}_2^2(\mathcal{L}(v_{m_n}),\mathcal{L}(v_{m_n}^{*} ))&=&{W}_2^2(\mathcal{L}(v_{m_n}-v),\mathcal{L}(v_{m_n}^{*}-v ))\\
\nonumber&\leqslant &\sup_{\varphi \in \mathcal{G} }\left\Vert \varphi^\prime_j\right \Vert_\infty^2 \frac 1 J \sum_{j=1}^J {W}_2^2(\mu_j,\mu_{n,j}).\end{eqnarray}
Now, recall that in the proof of Theorem \ref{th:testgen} one gets  that 
 ${W}_2^2(\mu_j,\mu_{n,j})= O_{\p} ( \frac{1}{\sqrt{n}} ) $. 
Then,
using that ${W}_r(\mathcal{L}(aX),\mathcal{L}(aY))=a{W}_r(\mathcal{L}(X),\mathcal{L}(Y))$ for $a>0$, \eqref{eq:cor1} gives
\begin{eqnarray} \label{eq:cor2}
\lefteqn{{W}_2^2\left(\mathcal{L}\left(\sqrt{m_n}\left( v_{m_n}-v\right)\right),\mathcal{L}\left(\sqrt{m_n}\left(v_{m_n}^{*}-v\right) \right)\right)}\hspace*{3cm}\\
\nonumber
&\leqslant& \frac{m_n}{ \sqrt{n}} \sup_{\varphi \in \mathcal{G} }\left\Vert  \varphi_j^\prime \right \Vert_\infty^2 \frac 1 J \sum_{j=1}^J  {n}{W}_2^2(\mu_j,\mu_{n,j}) \rightarrow 0
\end{eqnarray}
Moreover, under Assumptions   \ref{hyp:fct} to \ref{hyp:tassio12theta},   Theorem \ref{th:testgen}
gives  $\sqrt{m_n} \left(v^ 2_{m_n} - v^2 \right) \rightharpoonup \gamma$. 
If $v>0$, the classical Delta Method (see for instance in \cite{van2000asymptotic} p.25) gives 
$$\sqrt{m_n} \left(v_{m_n} - v \right) \rightharpoonup \frac{1}{2{v}}\gamma.$$
Hence \eqref{eq:cor2} enables to say that that $$\sqrt{m_n} \left(v^ *_{m_n} - v \right) \rightharpoonup \frac{1}{2{v}}\gamma.$$
Applying again a Delta Method leads to
$$\sqrt{m_n} \left((v^*)^2_{m_n} - v^2 \right)=\sqrt{m_n} \left(\inf_{\mathcal{G}}U^*_{m_n} - \inf_{\mathcal{G}}U \right) \rightharpoonup \gamma.$$
\eqref{eq:res_cor} is obtained by using Glivenko Cantelli Theorem and convergence of the empirical quantiles.

\hfill $\Box$

\bibliographystyle{imsart-nameyear}
\bibliography{biblithese}

\end{document}